\documentclass{amsart}
\usepackage[hypertex]{hyperref}
\makeatletter

\theoremstyle{plain}    
\newtheorem{thm}{Theorem}[section]
\numberwithin{equation}{section} 
\numberwithin{figure}{section} 
\theoremstyle{plain}    
\newtheorem{lem}[thm]{Lemma} 
\theoremstyle{plain}    
\newtheorem{prop}[thm]{Proposition} 
\theoremstyle{plain}    
\newtheorem{Def}[thm]{Definition} 
\theoremstyle{remark}
\newtheorem{rem}[thm]{Remark}
\theoremstyle{remark}

\usepackage[arrow,matrix,curve,ps]{xy}
\xyoption{dvips}
\CompileMatrices

\makeatother

\begin{document}

\pagestyle{myheadings}

\title{A Note on Invariants of Flows induced by Abelian Differentials on Riemann 
Surfaces}

\begin{abstract}
The real and imaginary part of any Abelian differential on a compact Riemann surface
define two flows on the underlying compact orientable
$\mathcal{C}^\infty$
surface. Furthermore, these flows induce an interval exchange transformation on every 
transversal simple closed curve, via Poincar\'e recurrence. This note shows that the 
ordered 
$K_0$-
groups of several 
$C^\ast$-
algebras naturally associated to one of the flows resp. interval exchange transformations
are isomorphic, mainly using the methods of I. Putnam~\cite{Put89,Put92}.
\end{abstract}


\author{Thomas Eckl}




\address{Thomas Eckl, Institut für Mathematik, Universität
  Bayreuth, 95440 Bayreuth, Germany}

\email{thomas.eckl@uni-bayreuth.de}

\urladdr{http://btm8x5.mat.uni-bayreuth.de/\~{}eckl}

\maketitle

\bibliographystyle{alpha}

\markboth{THOMAS ECKL}{INVARIANTS OF FLOWS INDUCED BY ABELIAN DIFFERENTIALS}

\setcounter{section}{-1}

\section{introduction}

\noindent
The real and imaginary part of an Abelian differential on a compact Riemann surface are
$\mathcal{C}^\infty$
real
$1$
forms defining flows with finitely many singularities of saddle type on the underlying
 compact orientable 
$\mathcal{C}^\infty$
surface
$M$.
Furthermore, each of these flows induces an interval exchange transformation on a closed
transversal curve, via Poincar\'{e} recurrence. If the Abelian differential is 
sufficiently general the induced flows are minimal. Then it is possible to associate
several 
$C^\ast$-
algebras to them: The crossed product
$C(M)\! >\!\!\!\triangleleft\ \mathbb{R}$
of the flow, the crossed product 
$C_0(M_0)\! >\!\!\!\triangleleft\ \mathbb{R}$
where
$M_0 = M - \{\mathrm{singularities}\}$
and the crossed product
$C(\Sigma)\! >\!\!\!\triangleleft\ \mathbb{Z}$
corresponding to a discrete dynamical system on a Cantor set
$\Sigma$
induced by the interval exchange transformation. Putnam~\cite{Put89,Put92} showed how to
compute the ordered 
$K_0$-\!
groups of
$C_0(M_0)\! >\!\!\!\triangleleft\ \mathbb{R}$
and 
$C(\Sigma)\! >\!\!\!\triangleleft\ \mathbb{Z}$
via embeddings in and inclusions of AF-algebras which induce order-isomorphies on the
$K_0$-
groups.

\noindent
The aim of this paper is to show the order-isomorphy of these 
$K_0$-
groups. Furthermore we prove that the cone defining the order is dual to the cone of 
invariant measures w.r.t. the interval exchange transformation. Both aims are achieved 
by using Veech's induced interval exchange transformations and Putnam's ``method of 
towers''. All these facts considered it is reasonable to expect the ordered groups to yield
interesting invariants of a Riemann surface together with an Abelian differential on 
it. It seems that Nikolaev \cite{Nik00, Nik01} tries to define and exploit such invariants but 
further clarification is certainly needed.

\noindent
To give a short synopsis of the article, the first two sections carefully explain the
necessary notions from the theory of dynamical systems and applies them on flows 
coming from Abelian differentials on Riemann surfaces, using classification results of 
Strebel \cite{Str84}. These sections contain more details than necessary for the 
specialists but try to be also readable for those mainly interested in new invariants of
Abelian differentials on Riemann surfaces. Following Veech \cite{Vee78,Vee82} 
section~\ref{IndIET-ssec} develops the notion of induced interval exchange 
transformations which is central to what follows. Section~\ref{crpr-sec} gives the 
details of Putnam's ``method of towers'' used to compute the ordered 
$K_0$-\!
groups of
$C_0(M_0)\! >\!\!\!\triangleleft\ \mathbb{R}$
and 
$C(\Sigma)\! >\!\!\!\triangleleft\ \mathbb{Z}$.
These details are needed in section~\ref{isoK0-sec} to produce a comprehensive proof of
the order-isomorphy of all occuring ordered groups and to show that the order-defining 
cone is dual to the cone of measures invariant w.r.t. the interval exchange 
transformation.

\section{Basics on flows}    \label{flow-sec} 

\begin{Def}
Let
$M$
be a compact orientable 
$\mathcal{C}^\infty$
surface. A 
$\mathcal{C}^\infty$
flow
$f^t$
is a
$\mathcal{C}^\infty$
map
$f : M \times \mathbb{R} \rightarrow M$
such that 
$f(m,0) = m$
and
$f(m,t_1+t_2) = f(f(m,t_1),t_2)$
for all
$m \in M$,
$t_i \in \mathbb{R}$.
\end{Def}

\noindent
In this note we only consider flows defined by closed real-valued
$1$-
forms 
$\phi_1$,
$\phi_2$
which are the real and imaginary part of a holomorphic
$1$-
form (Abelian differential)
$\omega = \phi_1 + i\phi_2$
on a Riemann surface
$C$.
These flows can be constructed as follows:
Around non-vanishing points of
$\omega$
it is trivial to find vector fields
$w_1$
such that
$\phi_1(w_1) = 0$
and
$\phi_2(w_1) > 0$
outside the zeroes, and around zeroes choose a holomorphic coordinate chart 
$V \rightarrow \mathbb{R}^2$
such that
$\omega = z^kdz$, 
and set
\[ w_1 = r \sin(k+1)\theta\frac{\partial}{\partial r} + 
         \cos(k+1)\theta\frac{\partial}{\partial \theta}. \]
Glue the finitely many vector fields by a partition of unity. In the same way it
is possible to construct a vector field 
$w_2$
on
$M$
such that
$\phi_1(w_2) > 0$
and
$\phi_2(w_1) = 0$.
The pair
$(w_1,w_2)$
is positively oriented outside the zeroes and in particular the
$w_i$
do not vanish outside the zeroes. Hence these vector fields can be integrated to
flows
$f_i^t$
having their fixed points (singularities) exactly in the zeroes of
$\omega$.

\noindent
Let
$f^t$
be a flow on a compact surface 
$M$.
The curve
$l(x) := \{f(x,t): -\infty < t < +\infty \}$
is called the \textit{trajectory} of
$f^t$
through
$x \in M$.
A trajectory
$l(m) = \{m\}$
is called a \textit{singular point}.

\noindent
There are two other types of trajectories, called \textit{regular}: if there 
exist
$t_1 \neq t_2$
such that 
$f(x,t_1) = f(x,t_2)$,
the trajectory through
$x$
is \textit{periodic} or \textit{closed}. Otherwise it is called 
\textit{non-closed}.

\noindent
For non-closed curves the rays
$l(x)^+ := \{f(x,t): 0 \leq t < \infty\}$
and
$l(x)^- := \{f(x,t): -\infty < t \leq 0\}$
are called the \textit{positive} resp. \textit{negative trajectory ray} through
$x$.
The
\textit{$\omega$-limit set} is defined as
\[ \omega[l(x)] := \{\tilde{x} \in M : \exists\ \mathrm{sequence\ } 
   (t_n)_{n \in \mathbb{N}},\ t_n \rightarrow \infty, 
   f(x,t_n) \rightarrow \tilde{x}\}, \]
and similarly the \textit{$\alpha$-limit set}
$\alpha[l(x)]$, 
with
$t_n \rightarrow -\infty$.

\noindent
A segment
$\Sigma: [0,1] \rightarrow M$
is called \textit{contact-free} or \textit{transversal} to the flow
$f^t$ 
iff for all
$m \in \Sigma \setminus \partial \Sigma$
there exists a neighborhood
$U$
of
$m$
and a diffeomorphism
$\phi: U \rightarrow \mathbb{R}^2$
such that
\[ \phi(f^t(m^\prime)) = \phi(m^\prime) + (t,0),\ \phi(m) = (0,0) \]
and
$\phi(\Sigma \cap U) = \{0\} \times [-1,1]$.
Since such a diffeomorphism exists around any non-singular (or 
\textit{regular}) point of a flow, we
have
\begin{lem} \label{transseg-lem}
Through each regular point of a flow there passes a contact-free segment.
\hfill $\Box$
\end{lem}

\noindent
A simple closed curve
$C$
(the image of an embedding of a circle in
$M$)
is called a \textit{contact-free cycle} or a \textit{closed transversal} of the
flow if its arcs are contact free segments. It is not hard to give an example 
of a flow for which no regular point has a contact free cycle passing through 
it (for example a flow on a sphere with two saddles and four centers). On the 
other hand there is an easy criterion for the existence of such a contact-free
cycle:
\begin{lem}[{\cite[Lem.2.1.2]{ABZ96}}] \label{transcyc-lem}
Suppose that the trajectory
$l$
of the flow
$f^t$
intersects a contact-free segment
$\Sigma$
at more than one point. Then there exists a contact-free cycle that intersects 
$l$.
\end{lem}

\noindent
The proof of this lemma uses the following 
\begin{thm}[Long Flow Tube Theorem]
Let
$d$
be a compact arc of a regular trajectory of a
$\mathcal{C}^\infty$-\!
flow
$f^t$,
and suppose that
$d$
does not form a closed curve. Then there exists a neighborhood 
$U$
of
$d$
and a
$\mathcal{C}^\infty$-\!
diffeomorphism
$\psi: U \rightarrow \mathbb{R}^2$
carrying the arcs in
$U$
of trajectories of
$f^t$
into trajectories of the dynamical system
$\dot{x} = 1, \dot{y} = 0$. 
\end{thm}

\noindent
For a proof see e.g. \cite{ALGM73}. 

\noindent
Let
$\overline{m_1 m_2}$
be the arc with endpoints
$m_1$
and
$m_2$
on a trajectory of the flow
$f^t$,
and let
$\Sigma_1$
and
$\Sigma_2$
be disjoint contact-free segments passing through
$m_1$
and
$m_2$,
respectively, such that
$\overline{m_1 m_2} \cap (\Sigma_1 \cup \Sigma_2) = \{m_1,m_2\}$.
For definiteness we will assume that 
$m_2 \in l(m_1)^+$.

\noindent
By the Long Flow tube Theorem, there exists a neighborhood
$\Sigma \subset \Sigma_1$
of
$m_1$
on
$\Sigma_1$
such that for any
$m \in \Sigma$
the positive semitrajectory
$l(m)^+$
intersects
$\Sigma_2$
without first intersecting
$\Sigma_1$.
Denote by
$\tilde{m}$
the first point where
$l(m)^+$
intersects
$\Sigma_2$.
\begin{Def}
The mapping
$P(m,\Sigma): \Sigma_1 \rightarrow \Sigma_2$
assigning the point 
$\tilde{m} \in \Sigma_2$
to a point
$m \in \Sigma$
according to the rule above is called the \textit{Poincar\'e mapping} or 
\textit{first return mapping} (induced by the flow
$f^t$).
\end{Def}

\noindent
The Long Flow Tube Theorem gives us
\begin{lem}[Poincar\'e mapping theorem]
Let
$P(m,\Sigma): \Sigma_1 \rightarrow \Sigma_2$
be the Poincar\'e mapping of the contact-free segment 
$\Sigma_1$
into the contact-free segment 
$\Sigma_2$
induced by a 
$\mathcal{C}^\infty$-\!
flow
$f^t$.
Then
$P$
is a
$\mathcal{C}^\infty$-\!
diffeomorphism of
$\Sigma_1$ 
onto its range.
\hfill $\Box$
\end{lem}

\noindent
The only flows we want to consider are those induced by (the imaginary part of) holomorphic
$1$-
forms. Their behaviour is in some sense more regular than the behaviour of arbitrary flows 
on surfaces, see for example in~\cite[III.2]{ABZ96}. The simplifications are described 
in~\cite{Str84}. 

\noindent
Let
$f^t_\phi$
be a flow on a compact surface 
$M$
induced by the holomorphic $1$-form
$\phi$. 
It can be shown that closed trajectories never come alone: the trajectories 
through points in a neighborhood of the closed trajectory 
$l(x)$
are also closed. One
may even find an open neighborhood of
$l(x)$
on
$M$
which is homeomorphic to a ring domain
$\{z: r < |z| < R \} \subset \mathbb{C}$
such that the circles around
$0$
are mapped onto the closed trajectories (this is a consequence of the fact, 
that the flow
$f^t_\phi$
preserves the volume given by
$|\phi|$). There is a maximal ring domain 
associated to 
$l(x)$
which is uniquely determined if
$M$
is not a torus. 
\begin{thm}[{\cite[Thm.9.4]{Str84}}]
Let 
$f^t_\phi$
be a flow on a compact Riemann surface 
$M$
induced by a holomorphic
$1$-
form
$\phi$.
Every closed trajectory
$l(x)$
of
$f^t$
is embedded in a maximal ring domain 
$R$ 
swept out by closed trajectories
$l(x^\prime)$, 
the \textit{ring domain} associated with
$l(x)$.
It is uniquely determined except for an orientable foliation with closed 
trajectories on a torus. The boundary of the ring domain consists of compact
critical leaves of the foliation.

\noindent
Two ring domains
$R_0$
and
$R_1$
associated with the closed trajectories
$l(x_0)$
and
$l(x_1)$
respectively are either disjoint or identical. If a closed trajectory
$l(x_1)$
is freely homotopic to
$l(x_0)$
then
$l(x_1) \subset R_0$,
which is equivalent to
$R_1 = R_0$.
\end{thm}

\noindent
Strebel~\cite[p.44pp]{Str84} gives a description of all possible types of 
non-closed trajectory rays and the structure of their limit sets: If
$\omega[l(x)] = \{P\}$,
then
$P$
is a singularity of the foliation (this is an easy exercise).
$l(x)^+$
is called a \textit{critical ray}, and there are only finitely many critical 
rays.

\noindent
If
$\omega[l(x)]$
contains more than one point, then it contains the initial point
$x$,
and
$\overline{l(x)} = \omega[l(x)]$
(this is again a consequence of the volume preserving property of the flow). 
Such a ray is called \textit{recurrent}. A trajectory both rays of which are
recurrent is called a \textit{spiral}.

\noindent
The interior of the limit set of a recurrent ray
$l(x)^+$
is non-empty and connected, hence a domain. Its boundary consists of compact 
critical rays of the flow and their limiting singularities. Every 
trajectory ray
$l(x^\prime)^+$
through a point
$x^\prime$
in the interior of
$\omega[l(x)]$
is everywhere dense in
$\omega[l(x)]$,
and its limit set coincides with
$\omega[l(x)]$.
If the limit sets of two recurrent rays have an interior point in common, they
coincide.

\noindent
Takink all these facts together we get
\begin{thm} \label{decompflow-thm}
Let 
$f^t$
be a flow on a compact Riemann surface 
$M$
induced by the imaginary part of a holomorphic
$1$-
form.
Then the connected components of
$M$
minus the compact critical rays of 
$f^t$
are bounded by critical segments, and either
\begin{itemize}
\item[(i)]
the components are ring domains swept out by closed trajectories, or
\item[(ii)]
they are the closure of a spiral trajectory.
\end{itemize}
\end{thm}

\noindent
The most general type of flow is given by
\begin{Def}
Let
$M$
be a compact orientable 
$\mathcal{C}^\infty$
surface. A 
$\mathcal{C}^\infty$
flow
$f^t$
is \textbf{minimal} iff
\begin{itemize}
\item[(i)]
every trajectory which is not a fixed point, is dense in
$M$, 
and
\item[(ii)]
the set of singularities 
$\mathrm{Sing} f^t$
consists of saddles with 
$2m$
separatrices.
\end{itemize}
\end{Def}

\noindent
This means in particular that every trajectory is non-closed and has a direction in which 
it is recurrent. There are only finitely many trajectories which are not spirals. 
Vice versa, it follows from theorem~\ref{decompflow-thm} that a flow on a compact Riemann 
surface induced by a holomorphic
$1$-
form without critical segments is minimal.

\section{Interval exchange transformations} \label{IET-sec}

\noindent
Let
$x$
be a regular point of a minimal flow
$f^t$
on a compact Riemann surface
$X$
of genus
$g \geq 2$
induced by the imaginary part of a holomorphic
$1$-
form on
$X$, 
and let
$\Sigma$
be a contact-free segment through 
$x$
(which exists by lemma~\ref{transseg-lem}). Since the trajectory 
$l$
through
$x$
is dense on
$X$,
the trajectory
$l$
intersects
$\Sigma$
in another point
$x^\prime$.
Consequently, by lemma~\ref{transcyc-lem} there exists a contact-free cycle
$C$.

\noindent
Let
$l(x)^+$
be a critical ray through a point
$x \in X$.
Since
$l(x)^-$
is dense in
$X$,
this ray intersects
$C$.
By a compactness argument, we can choose
$x^\prime \in l(x)$
such that
$l(x^\prime)^+ \cap C = \{x^\prime\}$.
Since there are only finitely many critical rays, we have
$n$
points
$x_1, \ldots, x_n \in C$
such that
$l(x_i)^+$
is critical and 
$l(x_i)^+ \cap C = \{x_i\}$.
Order
$x_1, \ldots, x_n \in C$
such that
$x_{i+1}$
is next to
$x_i$.

\noindent
Obviously, the Poincar\'e return mapping maps the open arc
$\overline{x_i x_{i+1}}$
to another arc
$\overline{y_i y_{i+1}}$
on
$C$.
The
$y_i$'s 
are points on critical rays of type
$l(y)^-$
such that
$l(y_i)^- \cap C = \{y_i\}$,
and the inverse mapping from
$\overline{y_i y_{i+1}}$
to
$\overline{x_i x_{i+1}}$
is the  
Poincar\'e return mapping of the reversed flow
$f^{-t}$.

\noindent 
Now, if
$\phi$
is the imaginary part of the holomorphic
$1$-
form defining
$f^t$
then the line element
$|\phi|$
defines a transverse measure to the trajectories of
$f^t$.
The Long Tube Theorem shows that the (transverse) length of the two arcs above 
remains equal. Cutting
$C$
at some point
$x_0$,
we get a map of the following type (at least outside the critical points 
$x_1, \ldots, x_n$): 
\begin{Def}
Let
$n \geq 2$
be an integer and
$\sigma \in \Sigma_n$
a permutation of
$\{1,2, \ldots,n\}$.
Let
$\alpha = (\alpha_1, \alpha_2, \ldots, \alpha_n)$
be an element of
$\mathbb{R}^n$.
For
$i = 0, \ldots, n$
define
\[ \beta_\alpha(i) = \sum_{j \leq i} \alpha_j,\   
    \beta^\prime_\alpha(i) = \sum_{\sigma(j) \leq i} \alpha_j. \]
Let 
$I_\alpha(i) = [\beta_\alpha(i-1), \beta_\alpha(i))$
and
$I^\prime_\alpha(i) = [\beta^\prime_\alpha(i-1),\beta^\prime_\alpha(i))$,
and let
\[ \tau_\alpha(i) = \sum_{\sigma(j) < \sigma(i)} \alpha_j - \sum_{j < i} \alpha_j, \]
for
$i = 1, \ldots, n$.
The bijective mapping
$T = T(\sigma, \alpha)$
of
$[0,\beta_\alpha(n))$
onto itself defined by
\[ T(x) = x + \tau_\alpha(i),\ x \in I_\alpha(i) \]
is called an \textbf{interval exchange transformation}, with data
$\sigma, \alpha$,
mapping
$I_\alpha(i)$
to
$I^\prime_\alpha(i)$.
\end{Def}

\noindent
Roughly speaking,
$T = T(\sigma, \alpha)$
partitions 
$[0,\beta_\alpha(n))$
into
$n$
subintervals of length
$\alpha_1, \ldots, \alpha_n$
and translates each so that they fit together to make up
$[0,\beta_\alpha(n))$
in a new order determined by
$\sigma$.

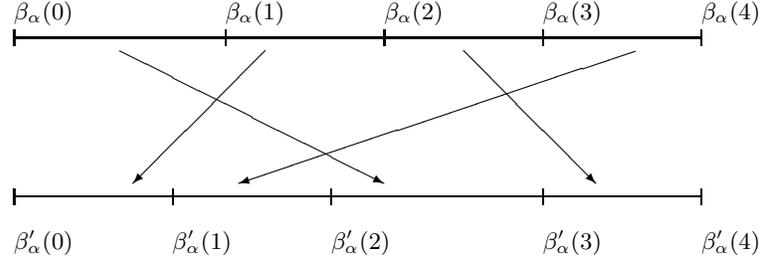
\begin{center}
\begin{picture}(270,100)(0,0)

\put(0,20){\line(1,0){260}}
\put(0,80){\line(1,0){260}}

\put(0,17){\line(0,1){6}}
\put(60,17){\line(0,1){6}}
\put(120,17){\line(0,1){6}}
\put(200,17){\line(0,1){6}}
\put(260,17){\line(0,1){6}}
\put(0,77){\line(0,1){6}}
\put(80,77){\line(0,1){6}}
\put(140,77){\line(0,1){6}}
\put(200,77){\line(0,1){6}}
\put(260,77){\line(0,1){6}}

\put(40,75){\vector(2,-1){100}}
\put(95,75){\vector(-1,-1){50}}
\put(170,75){\vector(1,-1){50}}
\put(235,75){\vector(-3,-1){150}}

\put(0,0){\small $\beta^\prime_\alpha(0)$}
\put(60,0){\small $\beta^\prime_\alpha(1)$}
\put(120,0){\small $\beta^\prime_\alpha(2)$}
\put(200,0){\small $\beta^\prime_\alpha(3)$}
\put(260,0){\small $\beta^\prime_\alpha(4)$}

\put(0,87){\small $\beta_\alpha(0)$}
\put(80,87){\small $\beta_\alpha(1)$}
\put(140,87){\small $\beta_\alpha(2)$}
\put(200,87){\small $\beta_\alpha(3)$}
\put(260,87){\small $\beta_\alpha(4)$}

\end{picture}
\end{center}

\noindent
The minimality of the flow
$f^t$
with which we started above implies the minimality of
$T$,
that is, the orbit of every point in
$[0,\beta_\alpha(n))$
is dense. (The orbit of a point
$x$
under
$T$,
or the
\textit{$T$-\!
orbit} of
$x$,
is
$\{T^n(x)\ |\ n \in \mathbb{Z}\}$.)
A necessary condition for minimality is that there is no
$k$
such that
$\sigma(\{1, \ldots, k\}) = \{1, \ldots, k\}$. 
In this case,
$\sigma$
is called \textit{irreducible}.
A (non-trivial) sufficient condition (\cite{Kea75}) is given by 
the following
\begin{Def} \label{IDOC-def}
$T = T(\sigma, \alpha)$
satisfies the infinite distinct order condition (for short IDOC) iff 
$\sigma$
is irreducible and the 
$T$-\!
orbits of the points
$\beta_\alpha(1), \beta_\alpha(2), \ldots, \beta_\alpha(n-1)$
are all infinite and distinct. 
\end{Def}

\noindent
The interval exchange transformation which we constructed from (the imaginary
part of) a holomorphic
$1$-
form 
$\omega$
on
$X$
and a closed transversal 
$C$
still depends on the point where we cut 
$C$.
If we choose one of the 
$y_i$
(see above) we have (notations as above)
\begin{lem} \label{FolIET1-lem}
Let
$\omega$
have 
$l$
zeroes
$p_1, \ldots, p_l \in X$,
and let
$T = T(\sigma, \alpha)$
be the induced interval exchange transformation when cutting 
$C$
in one of the points
$y_i$. 
Then 
$T$
exchanges 
$2g + l - 1$
intervals.
\end{lem}
\begin{proof}
Let
$k_1, \ldots, k_l$
be the multiplicities of the zeroes
$p_1, \ldots, p_l$.
Then the singularity 
$p_i$
of the flow induced by 
$\mathrm{Im}\ \omega$
is a
($2k_i + 2$)-\!
prong saddle, hence there are
$k_i + 1$
critical rays emanating from
$p_i$. 
The construction of
$T$
implies that there are
$k_1 + \cdots + k_l + l$
points
$x_i$,
hence intervals on
$C$.
Cutting in
$y_i$
divides one of these intervals in two parts since no
$y_i$
can be equal to some
$x_j$
(otherwise, there exists a critical segment, and the foliation is not minimal). On 
the other hand, no interval is mapped to an interval containing
$y_i$
in its interior. Consequently, 
$T$
exchanges
$k_1 + \cdots + k_l + l + 1$
intervals.

\noindent
Finally
$k_1 + \cdots + k_l = 2g - 2$
by Riemann-Roch, and the lemma follows.
\end{proof} 

\noindent
That
$T$
comes from cutting the closed curve 
$C$
is reflected by the condition
\[ \sigma(1) = j \Rightarrow \sigma(n) = j-1. \]
Here,
$j$
is always
$\neq 1$
because otherwise
$\sigma$
would be reducible. A
$T$
constructed as above is called the intervall exchange transformation \textit{induced}
by the holomorphic
$1$-
form
$\omega$
(and the simple closed curve
$C$).
It satisfies the IDOC because otherwise there would exist a critical segment.

\noindent
Putnam (\cite{Put92}) described a procedure how to construct a flow 
$F$
on a compact oriented surface 
$M$
resolving the discontinuities of an interval exchange transformation
$T$,
together with a closed transversal
$N$ 
such that the Poincar\'e return mapping on
$N$
is just the interval exchange
$T$,
except at the discontinuities of
$T$
which flow directly into the singularities:

\noindent
Begin with 
$P = [0,1) \times [0,1)$.
Let
$V^\prime(0) = V(0) = \{0\} \times [0,1]$
and
$V^\prime(n) = V(n) = \{1\} \times [0,1]$.
For
$i = 1, \ldots, n-1$,
let
$V(i) = (\beta(i),1)$
and
$V^\prime(i) = (\beta^\prime(i),0)$.
We define
$M$
to be the quotient of
$P$
obtained by collapsing
$V(0)$
and
$V(n)$
to single points and by identifying
$I(i) \times \{1\}$
with
$I^\prime(i) \times \{0\}$
via
$T$, 
for
$i = 1, 2, \ldots, n$.

\noindent
Define
$\sigma_0$,
a permutation of
$\{0,1,\ldots,n\}$
by
\[ \sigma_0(j) = \left\{ \begin{array}{cl}
   \sigma^{-1}(1)-1 & \mathrm{if\ } j=0, \\
   n                & \mathrm{if\ } j=\sigma^{-1}(n), \\
   \sigma^{-1}(\sigma(j)+1)-1 & \mathrm{otherwise.}
                         \end{array} \right.\]
We let
$N(\sigma)$
denote the number of cyclic components of
$\sigma_0$
and we let
$X(\sigma)$
denote the image of
$\{V(0), \ldots, V(n), V^\prime(0), \ldots, V^\prime(n)\}$
in
$M$.
Then it is easy to check that if
$\sim$
denotes the equivalence relation generated by
$V(i) \sim V(\sigma_0(i))$,
$i = 0, \ldots, n$,
then 
$V(i) \sim V(j)$
iff they have the same image in
$M$.
It is clear that each
$V^\prime(i)$
has the same image as some
$V(j)$. 
Thus we may identify
$X(\sigma)$
with
$\{V(0), \ldots, V(n)\}/\sim$
and
$X(\sigma)$
just consists of
$N(\sigma)$
points in
$M$.

\noindent
The flow
$F: \mathbb{R} \times M \rightarrow M$
is obtained by integrating the vector field
$X(p) = \omega(p) \cdot (0,1)$,
for
$p \in M$,
where
$\omega(p)$
is a certain positive scalar function on
$M$
vanishing exactly on
$X(\sigma)$
and
$(0,1)$
comes from the obvious vector field on
$P$.
A transversal 
$N$
is given by the image in
$M$
of
$[0,1] \times \{\frac{1}{2}\}$.
It is closed iff 
$V(0)$
is identified with
$V(n)$. 
This is the case if
$\sigma(1) = j$ 
and
$\sigma(n) = j-1$
for some
$j > 1$,
since then
$\sigma_0(n) = 0$.

\begin{lem} \label{1-form-IET-lem}
Let
$X$
be a Riemann surface of genus
$g \geq 2$,
let
$C$
be a simple closed curve transversal to the minimal flow
$f^t$
induced by the imaginary part of a holomorphic 
$1$-\!
form
$\omega$
on
$X$.
Let
$T(\sigma,\alpha)$
be the induced interval exchange transformation, and let
$F^t$
be the flow constructed on
$M$
from
$T(\sigma,\alpha)$
as above. Then there exists a homeomorphism
$\phi: X \rightarrow M$
equivariant with respect to the flows 
$f^t$
and
$F^t$.
In particular, the 
$X(\sigma)$
points on
$M$
correspond to the zeroes of
$\omega$,
and the length of the cycles belonging to one of these
$N(\sigma)$
points (after possibly omitting 
$0$
and
$n$) 
equals the multiplicity
$+ 1$
of the corresponding zeroes. The closed transversal
$N$
is homotopic to
$C$. 
\end{lem}
\begin{proof}
We easily check that the singularities of
$F$
on
$M$
are saddles, where the emanating rays are counted by the members of the cycle belonging
to the singularity (they consist of the points 
$V(i)$
which are identified with the singularity on 
$M$)
and the incoming rays by the points
$V^\prime(j)$
which are identified with the singularity on
$M$. 
The only exceptions are 
$V(0) = V^\prime(0)$
and
$V(n) = V^\prime(n)$
from which no ray emanates resp. to which no ray in 
$P$
runs in.

\noindent
Now look at the singularity
$p$
to which
$V(0)$
and
$V(n)$
are mapped.
$N$ 
runs through
$p$
but if we perturb
$N$
appropriately in a small neighborhood around
$p$, 
we get a closed \textit{everywhere} transversal curve
$N^\prime$.
The construction of 
$T(\sigma,\alpha)$
from
$\omega$
and
$C$
(and especially lemma~\ref{FolIET1-lem}) shows that we can identify
$C$
and
$N^\prime$
and the intervalls given on them. This identification may be expanded to the wished 
equivariant homeomorphism
$\phi: X \rightarrow M$.
\end{proof}

\noindent
This lemma allows to exclude that 
$T$
may be represented as a transformation exchanging less intervals:
\begin{lem}
Let
$T = T(\sigma,\alpha)$
be the the interval exchange transformation \textit{induced}
by the holomorphic
$1$-
form
$\omega$
and the simple closed curve
$C$. 
Then there is no
$1 \leq j < n$
such that
\[ \sigma(j+1) = \sigma(j) + 1. \]
\end{lem}
\begin{proof}
Obviously,
\[ \sigma_0(j) = \sigma^{-1}(\sigma(j)+1) - 1 = \sigma^{-1}(\sigma(j+1)) - 1 = j. \]
This means that 
$V(j)$
would be a 2-prong saddle of the flow 
$f^t$,
and cannot be a singularity of (the imaginary part of)
$\omega$:
contradiction.
\end{proof}

\section{Induced interval exchange transformations} \label{IndIET-ssec}

\noindent
This construction was introduced by Keane~\cite{Kea75} and
systematically used by Veech~\cite{Vee78, Vee82} to prove ergodicity results for interval
exchange transformations. It starts with the following
\begin{Def}
An interval
$J = [a,b) \subset [0, \beta(n))$
is \textit{admissible} for
$T = T(\sigma,\alpha)$
if 
$a = T^k\beta_\alpha(i)$,
$b = T^l\beta_\alpha(j)$
for some
$1 \leq i,j < n$
and 
$k,l$
satisfying either
\begin{itemize}
\item[(a)]
$k \geq 0$, 
and there is no
$m$
such that 
$0 < m < l$
and
$T^m\beta_\alpha(i) \in J$
(and similarly for
$l$),
or
\item[(b)]
$k < 0$,
and there is no
$m$,
$k < m \leq 0$,
such that
$T^m\beta_\alpha(i) \in J$
(and similarly for
$l$).
\end{itemize}
\end{Def}

\noindent
In particular, every interval
$[\beta_\alpha(i),\beta_\alpha(i+1))$
is admissible (here,
$k = l = 0$).

\noindent
Assume now that
$T = T(\sigma,\alpha)$
satisfies the IDOC (see definition~\ref{IDOC-def}), and let
$J = [a,b)$
be an admissible interval for
$T$. 
We define 
$U$
to be the induced first return transformation of
$T$
on
$J$.
For each
$j$
such that 
$1 \leq j \leq n$, 
the minimality of
$T$
implies that there exists a least
$k_j \geq 0$
such that
$T^{-k_j}\beta_\alpha(j) \in J$.
This determines
$m - 1$
distinct points by the IDOC, and the distances between consecutive points from left to 
right (including 
$a$
and
$b$)
determine
$\alpha^\prime = (\alpha^\prime_1, \ldots, \alpha^\prime_n) \in \mathbb{R}^n$ 
such that
$[a,a + \beta_{\alpha^\prime}(n)) = J$
and
$T^{-k_{j_0}}\beta_\alpha(j_0) = \beta_{\alpha^\prime}(i_0)$,
where
$i_0$
denotes the position of
$T^{-k_{j_0}}\beta_\alpha(j_0)$
in the
$m-1$
distinct points
$T^{-k_j}\beta_\alpha(j) \in J$.

\noindent
It is also true for
$1 \leq j \leq n$
that
$T^m (a + I_{\alpha^\prime}(j)) = T^m [a + \beta_{\alpha^\prime}(j-1),
                                       a + \beta_{\alpha^\prime}(j))$
lies in some
$I_\alpha(i)$
for
$0 \leq m < \mathrm{time\ of\ return\ to\ } J$.
For given
$j,i$
we denote by
$A_{ij}$
the number of times 
$T^m (a + I_{\alpha^\prime}(j)) \subset I_\alpha(i)$
for
$0 \leq m < \mathrm{time\ of\ return\ to\ } J$.
\begin{prop} \label{indIET-prop}
Let
$T = T(\sigma,\alpha)$
be an interval exchange transformation satisfying the IDOC, let
$J = [a,b) \subset [0,\beta_\alpha(n))$
be an admissible interval, and let
$U$,
$\alpha^\prime$
and
$A = (A_{ij}) \in M(n,\mathbb{Z})$
be as above. Then there exists a permutation
$\sigma^\prime \in \Sigma_n$
such that
$U - a = T^\prime = T(\sigma^\prime, \alpha^\prime)$
is an interval exchange transformation on
$[0, b-a)$
satisfying the IDOC. Furthermore
\[ \det A = \pm 1,\ \ \ \alpha = A \alpha^\prime. \]
\end{prop}
\begin{proof}
We can repeat all the arguments in~\cite[\S 3]{Vee78}. The only difference is that the 
separation points
$\beta_\alpha$
of
$T$
need not be points of discontinuity, that is there may exist
$j \in \{1, \ldots, n\}$
such that
$\sigma(j+1) = \sigma(j) + 1$
(and similarly for
$\sigma^\prime$).

\noindent
The idea behind the equality
$\alpha = A \alpha^\prime$
is that
$I_\alpha(i)$
is composed of
$A_{ij}$
copies of the interval
$I_{\alpha^\prime}(j)$,
$1 \leq j \leq n$. 

\noindent
That 
$U$
satisfies the IDOC is a consequence of the IDOC for
$T$,
the representation of the separation points of
$U$
by iterated
$T$-
images of separation points of
$T$
and the representation of
$U$
restricted to the intervals by iterations of
$T$.
\end{proof}

\noindent
Induced interval exchange transformations are so important because they can be used to 
compute the (nonnegative 
nonzero) invariant measures for the interval exchange transformation
$T(\sigma, \alpha)$,
in particular the ergodic ones.
\begin{Def}
A nonnegative nonzero
$T$-
invariant Borel measure 
$\mu$
on
$[0, \beta(n))$
is called \textbf{ergodic} if 
$\mu(A) = 0$
or
$\mu(A) = \mu([0, \beta(n)))$
for any 
$T$-
invariant subset
$A \in [0, \beta(n))$.
\end{Def} 

\noindent
The set of all 
$T$-
invariant Borel measures on
$[0, \beta(n))$
is a cone
$\Sigma(\sigma, \alpha)$.
It contains the Lebesgue measure, and if
$T$
is minimal the Lebesgue measure is necessarily an ergodic invariant measure.
$T$
is called \textbf{uniquely ergodic} if every 
$T$-
invariant measure is a multiple of Lebesgue measure. 

\noindent
In general, it is well known that the extremals of
$\Sigma(\sigma, \alpha)$
are the ergodic invariant measures (see e.g. \cite[Thm.6.10]{Wal82}). Furthermore, a 
minimal transformation exchanging
$n$
intervals can have at most
$n$
different ergodic measures (\cite[Thm.5.2.1]{CFS82}). As
$\Sigma(\sigma, \alpha)$
is spanned by its extremals,
\[ \nu(\sigma, \alpha) = \dim \Sigma(\sigma, \alpha) \]
is the number of ergodic invariant measures  for
$T$.

\noindent
In the beginning it was conjectured that every minimal interval exchange transformation is
uniquely ergodic. But soon counterexamples were discovered (see \cite[\S 5.4]{CFS82}), and 
Keane conjectured instead that \textit{almost all} interval exchange transformations (with
respect to the Lebesgue measure on the vectors
$\sigma \in \mathbb{R}^n$)
are uniquely ergodic. This conjecture was proven simultaneously by Masur \cite{Mas82} and 
Veech \cite{Vee82}. Veech's idea was to use the induced interval exchange transformations 
to describe
$\Sigma(\sigma, \alpha)$:
\begin{prop} \label{ergcone-prop}
Let
$T = T(\sigma, \alpha)$
be a minimal interval exchange transformation, and let
$\Lambda_n \subset \mathbb{R}^n$
be the simplex generated by the unit vectors in
$\mathbb{R}^n$.
If
$J_1 \supset J_2 \supset \cdots$
is a sequence of admissible intervals shrinking to a point (or
$\emptyset$)
then
\[ \Sigma(\sigma, \alpha) \cong \bigcap_{i=1}^\infty (A_1 A_2 \cdots A_i \Lambda_n) \]
where the
$A_i$ 
are the matrices describing the transition from the
$(i-1)$st
to the
$i$th 
induced interval exchange transformation. In particular, the cone on the right is spanned 
by the set of cluster points of the images of the extremal rays of
$\Lambda_n$
under the sequence
\[ A^{(i)} = A_1 A_2 \cdots A_i. \] 
\end{prop}

\begin{rem}
A sequence of admissible intervals
$J_1 \supset J_2 \supset \cdots$
shrinks to 
$\emptyset$
iff the half open intervals 
$J_n$
have the same right end 
$b$
(for all
$n \geq N$). 
Of course it is reasonable to say that such a sequence shrinks to
$b$
(from the left). This point of view agrees with the splitting of orbit points
$T^k\beta(i)$
in the Cantor set interpretation of an interval exchange transformation (see 
section~\ref{crpr-sec}) because
$b$
is of this form.
\end{rem}

\begin{proof}
See \cite[Prop.3.22]{Vee78}. The isomorphism is given by the map
$\Sigma(\sigma, \alpha) \rightarrow \Lambda_n$
mapping a
$T$-
invariant measure
$\mu$
to 
\[ \lambda(\mu) = (\mu([0,\beta(1))), \mu([\beta(1), \beta(2))), \ldots, 
                                     \mu([\beta(n-1), \beta(n)))) \in \mathbb{R}^n. \]

\noindent
For the last statement we have to prove that
\[ A^{(i+1)} \Lambda_n \subset A^{(i)} \Lambda_n. \]
This is a direct consequence of the nonnegativity of the entries of
$A_i$
which implies that the 
$(i+1)$st
images of the extremal rays are positive linear combinations of the
$i$th
images.

\noindent
The inclusion of the cone on the right side in the image of this map follows from the fact
that each column of
$A^{(i)}$,
when normalized, gives the relative frequency of visits of a specific maximal interval in
$J_n$
to the intervals
$[\beta(k), \beta(k+1))$
before returning to
$J_n$. 
Therefore, any cluster point of the sequence of normalized columns corresponds to an 
invariant measure for
$T(\sigma, \alpha)$,
by a construction described in the next lemma. 
\end{proof}

\noindent
\begin{lem} \label{clustinv-lem}
Let
$T(\sigma, \alpha)$
be a minimal interval exchange transformation, let
$p \in [0, \beta(n))$
be any point. Take any sequence
$m_1, m_2, \ldots$
of integers
$\geq 0$
(which may be equal or not, increasing or not)
and an increasing sequence
$n_1 < n_2 < \ldots < n_k < \ldots$
of positive integers such that for each 
$1 \leq j \leq n$
the limit
\[ \sigma_j^\prime := \lim_{k \rightarrow \infty} 
                      \frac{1}{n_k} \sum_{i=m_k}^{m_k+n_k} \chi_{I_j}(T^i p) \]
exists. Then there exists a unique 
$T$-
invariant measure
$\mu$
on
$[0, \beta(n))$
such that
\[ \mu(I_j) = \sigma_j^\prime. \]
\end{lem}
\begin{proof}
Let 
$\delta_{T^i p}$
be the Dirac distribution with center
$T^i p$,
and set
\[ \mu_k := \frac{1}{n_k} \sum_{i=m_k}^{m_k+n_k} \delta_{T^i p}. \]
We know that 
$\lim_{k \rightarrow \infty} \mu_k(I_j) = \sigma_j^\prime$
and we want to show that the limit is a
$T$-
invariant measure. For the
$T$-
invariance we compute for any Borel set
$I \subset [0, \beta(n))$
such that
$\lim_{k \rightarrow \infty} \mu_k(I)$
exists
\[ \begin{array}{rcl}
   \lim\limits_{k \rightarrow \infty} \mu_k(TI) & = & 
   \lim\limits_{k \rightarrow \infty} \frac{1}{n_k} 
   \sum\limits_{i=m_k}^{m_k+n_k} \delta_{T^i p}(TI) =
   \lim\limits_{k \rightarrow \infty} \frac{1}{n_k} 
   \sum\limits_{i=m_k}^{m_k+n_k} \delta_{T^{i-1} p}(I) = \\
    & & \\
    & = & \lim\limits_{k \rightarrow \infty} \frac{1}{n_k} 
          (\delta_{T^{m_k-1} p}(I) - \delta_{T^{m_k+n_k} p}(I)) +
   \lim\limits_{k \rightarrow \infty} \frac{1}{n_k} 
   \sum\limits_{m_k=0}^{m_k+n_k} \delta_{T^i p}(I) = \\
    & & \\
    & = &  \lim\limits_{k \rightarrow \infty} \frac{1}{n_k} 
   \sum\limits_{m_k=0}^{m_k+n_k} \delta_{T^i p}(I) = 
   \lim\limits_{k \rightarrow \infty} \mu_k(I).
   \end{array} \]
Hence the limit exists for all 
$T^k I_j$.

\noindent
Now we can use that the
$T^k I_j$
form a Dynkin system for the Borel algebra of
$[0, \beta(n))$:
The limit exists for all intervals in the common refinements of the partitions 
$\{I_j\}, \ldots, \{T^k I_j\}$, 
and these intervals generate the Borel algebra. Thus 
$\mu := \lim_{k \rightarrow \infty} \mu_k$
exists. At the same time this shows the uniqueness.
\end{proof}

\noindent
The description of
$\Sigma(\sigma, \alpha)$
in proposition~\ref{ergcone-prop} led to an easy criterion for unique ergodicity:
\begin{prop}
Let
$T = T(\sigma, \alpha)$,
$J_1 \supset J_2 \supset \cdots$
and
$A_1, A_2 \ldots $ 
be as before. If there exists a matrix 
$B$
with positive entries such that for infinitely many
$i,j$
\[ A_i A_{i+1} \cdots A_j = B, \]
then
$T$
is uniquely ergodic.
\end{prop}
\begin{proof}
See \cite[Prop.3.30]{Vee78}.
\end{proof}

\noindent
The heart of Veech's proof in \cite{Vee82} was to construct an absolute continuous measure on
the space of all interval exchange transformations whose permutation are belonging to the
same so-called ``Rauzy class'' and which is invariant with respect to the map given by 
induction of interval exchange transformations.

\section{The crossed product $C^\ast-\!$algebras and their $K_0$ groups}
\label{crpr-sec}

\noindent
In this section we associate several
$C^\ast-\!$
algebras to a minimal flow 
$f^t$
on a compact oriented surface
$M$
induced by a holomorphic 
$1$-
form and to the induced interval exchange transformation
$T$.

\noindent
The most natural 
$C^\ast-\!$
algebra associated to
$f^t$
is the crossed product 
$C^\ast-\!$
algebra
\[ C(M) >\!\!\!\triangleleft_{f^t}\ \mathbb{R} \]
associated to the action of
$\mathbb{R}$
induced on
$M$
by
$f^t$
(for details about the crossed product construction see \cite[12.1]{GVF-NC}).

\noindent
Let
$f^t$
have the singularities in
$x_1, \ldots, x_n$
and set
$M_0 := M - \{x_1, \ldots, x_n\}$.
Since
$M_0$
is
$f^t$-
invariant we can consider the crossed product
$C^\ast$-
algebra
\[ C_0(M_0) >\!\!\!\triangleleft_{f^t}\ \mathbb{R}. \]
We have the following
$f^t$-
invariant short exact sequence
\[ 0 \rightarrow C_0(M_0) \rightarrow C(M) \rightarrow 
   C(\{x_1, \ldots, x_n\}) \rightarrow 0, \]
and taking crossed products we obtain the short exact sequence
\[ 0 \rightarrow C_0(M_0) >\!\!\!\triangleleft_{f^t}\ \mathbb{R} \rightarrow 
                 C(M) >\!\!\!\triangleleft_{f^t}\ \mathbb{R} \rightarrow 
                 C(\{x_1, \ldots, x_n\}) >\!\!\!\triangleleft_{f^t}\ \mathbb{R} 
     \rightarrow 0. \]
Here we know that 
\[ C(\{x_1, \ldots, x_n\}) >\!\!\!\triangleleft_F\ \mathbb{R} \cong 
   \bigoplus_{i=1}^l C_0(\mathbb{R}) \]
because the crossed product
$C_0(\{x\}) >\!\!\!\triangleleft\ \mathbb{R} \cong C_0(\mathbb{R})$.

\noindent
Next, let 
$T(\sigma,\alpha)$
be the interval exchange transformation induced on a closed simple curve by
$f^t$. 
Putnam (\cite{Put89}) attached a 
$C^\ast$-
algebra to 
$T$ 
by constructing a Cantor set
$\Sigma$
and a homeomorphism
$\phi$
of
$\Sigma$
so that
$[0,\beta(n))$
is densely contained in
$\Sigma$
(in a natural way) and
$\phi_{|[0,\beta(n))} = T$,
proceeding as follows: Let
$D(T)$
denote the 
$T$-
orbits of
$\beta(1), \ldots, \beta(n-1)$, 
omitting the point
$0$.
We want to consider the set
$D(T) \times \{0,1\}$, 
but it will be more convenient to denote
$(x,0)$
and
$(x,1)$
by
$x^+$
and
$x^-$,
respectively. Let
\[ \Sigma = [0,1] - D(T) \cup \{x^+,x^- | x \in D(T)\}. \]
(This amounts to inserting the ``Cantor gaps'' at the points of 
$D(T)$.) 
There is an obvious linear order on
$\Sigma$,
using
$x^- < x^+$,
for all
$x \in D(T)$.
Endowed with the order topology,
$\Sigma$
is a Cantor set since the minimality of
$T$
insures that
$D(T)$ 
is dense. We include
$[0,\beta(n))$
in
$\Sigma$
by mapping
$x$
in
$D(T)$
to
$x^+$.
The definition of 
$\phi$
and the fact that it is a homeomorphism are both clear. Defining an automorphism on 
$C(\Sigma)$
via
$f \mapsto f \circ \phi$
gives the crossed product
\[ C(\Sigma) >\!\!\!\triangleleft_\phi\ \mathbb{Z}. \] 

\noindent
This
$C^\ast$-
algebra may also be defined as an operator algebra on
$L^2([0,\beta(n)))$: 
The Lebesgue measure on
$[0,\beta(n))$
induces a 
$\phi$-
invariant Borel measure 
$\mu$
on
$\Sigma$
with respect to the order topology. Using the reduced crossed product description (see
\cite[Ex.12.1]{GVF-NC}) we conclude that the
$C^\ast$-
algebra of operators on 
$L^2(\Sigma,d\mu)$
generated by the multiplication operators
$\xi(x) \mapsto f(x)\xi(x)$
and the shift operator
$V\xi(x) := \xi(\phi^{-1}(x))$
is isomorphic to
$C(\Sigma) >\!\!\!\triangleleft_\phi\ \mathbb{Z}$. 
Since
$\phi$
is minimal, we can restrict the generating operators to
$V$
and the multiplication operators defined by the characteristic functions
$\chi_{[\beta(i),\beta(i+1)]}$,
$i = 0, \ldots ,n-1$
of the intervals involved in 
$T(\sigma,\alpha)$.
These operators can be transferred to operators on
$L^2([0,\beta(n)))$.

\noindent
In order to investigate its
$K$-
theory Putnam (\cite{Put89}) constructed an AF-subalgebra of 
$C(\Sigma) >\!\!\!\triangleleft_\phi\ \mathbb{Z}$
for every non-empty closed subset
$Y$
of
$\Sigma$
and proved (for definitions and first properties of AF-algebras see \cite{Dav96}) 
\begin{thm}[\cite{Put89}] \label{PutK0-thm}
Let
$Y$
be a non-empty closed subset of
$\Sigma$.
Then
$A_Y$, 
the
$C^\ast$-
subalgebra of 
$C(\Sigma) >\!\!\!\triangleleft_\phi\ \mathbb{Z}$
generated by 
$C(\Sigma)$
and
$VC_0(\Sigma-Y)$,
is an AF-algebra. If 
$i$
denote the inclusion map of
$A_Y$
in
$C(\Sigma) >\!\!\!\triangleleft_\phi\ \mathbb{Z}$
there is an exact sequence 
\[ 0 \rightarrow \mathbb{Z} \stackrel{\alpha}{\rightarrow} C(Y,\mathbb{Z}) \rightarrow 
   K_0(A_Y) \stackrel{i_\ast}{\rightarrow} 
   K_0(C(\Sigma) >\!\!\!\triangleleft_\phi\ \mathbb{Z}) \rightarrow 0 \]
where
$\alpha$
is the map taking
$n \in \mathbb{Z}$
to the constant function
$n$.

\noindent
Moreover, for every 
$a \in K_0^+(C(\Sigma) >\!\!\!\triangleleft_\phi\ \mathbb{Z})$,
there is
$b \in K_0^+(A_Y)$
such that
$i_\ast(b) = a$.

\noindent
In particular, if
$Y$
is a single point,
$i_\ast$
is an isomorphism of ordered groups.
\end{thm}

\noindent
Since we need the techniques used in the proof later on we reproduce them in all details:
Putnam explicitely constructed the Bratteli diagram of the 
AF-algebra
$A_Y$
(and hence its
$K_0$
group) by the ``method of towers''. To understand it, we first need some notation: We 
say that a subset
$E$
of
$\Sigma$
is clopen if it is both closed and open. We let
$\chi_E$
denote the characteristic function of
$E$,
which will be continuous if
$E$
is clopen. A partition
$\mathcal{P}$
of
$\Sigma$
is defined to be a finite collection of pairwise disjoint clopen sets whose union is all
of
$\Sigma$.
If
$\mathcal{P}$
is a partition of
$\Sigma$,
we let
$\mathcal{C}(\mathcal{P}) = \mathrm{span}\{\chi_E | E \in \mathcal{P}\}$.
$\mathcal{C}(\mathcal{P})$
may be viewed as those functions in
$C(\Sigma)$
which are constant on each element of
$\mathcal{P}$. 
The fact that
$\Sigma$
is totally disconnected implies that any function in
$C(\Sigma)$
may be approximated arbitrarily closely by one in some 
$\mathcal{C}(\mathcal{P})$.
Given two partitions
$\mathcal{P}_1$
and
$\mathcal{P}_2$
of
$\Sigma$,
we say 
$\mathcal{P}_2$
is finer than
$\mathcal{P}_1$
and write
$\mathcal{P}_2 \geq \mathcal{P}_1$,
if each element of
$\mathcal{P}_2$
is contained in a single element of
$\mathcal{P}_1$.
This is clearly equivalent to the condition that
$\mathcal{C}(\mathcal{P}_1) \subset \mathcal{C}(\mathcal{P}_2)$.
Given two partitions
$\mathcal{P}_1$
and
$\mathcal{P}_2$,
we define the partition 
$\mathcal{P}_1 \vee \mathcal{P}_2$
to be
$\{E \cap F | E \in \mathcal{P}_1, F \in \mathcal{P}_2\}$.

\noindent
The first step is to show that a partition 
$\mathcal{P}$
of
$\Sigma$
and a non-empty closed subset
$Y$
of
$\Sigma$
give rise to a finite-dimensional 
$C^\ast$-
subalgebra of
$C(\Sigma) >\!\!\!\triangleleft_\phi\ \mathbb{Z}$.
\begin{lem}
The
$C^\ast$-
subalgebra of
$C(\Sigma) >\!\!\!\triangleleft_\phi\ \mathbb{Z}$
generated by 
$\mathcal{C}(\mathcal{P})$
and
$V\chi_{\Sigma - Y}$
is finite dimensional.
\end{lem}
\begin{proof}
We begin by defining
$\lambda: Y \rightarrow \mathbb{Z}$
by
\[ \lambda(y) = \inf \{n \geq 1 | \phi^n(y) \in Y\},\ \ y \in Y. \]
Notice that since 
$\phi$
is minimal and
$Y$
is open, there is, for each point
$y$,
a positive integer
$n$
such that
$\phi^n(y) \in Y$,
so
$\lambda$
is well-defined. 

\noindent
It is straightforward to verify that
$\lambda$
is upper (lower) semi-continuous because
$Y$
is open (closed), and so
$\lambda$
is continuous. Then because
$Y$
is compact,
$\lambda(Y)$
is finite. Let us suppose that
$\lambda(Y) = \{J_1, J_2, \ldots, J_K \}$
with
$J_1 < \cdots < J_K$.

\noindent
For
$k = 1, \ldots, K$
and
$j = 1, \ldots, J_k$
define the clopen set
$Y(k,j) = \phi^j(\lambda^{-1}(J_k))$. 
It follows at once from the definitions that the following properties hold:
\begin{itemize}
\item[(1)]
$\bigcup_{k=1}^{K} Y(k,1) = \phi(Y)$,
\item[(2)]
$\phi(Y(k,j)) = Y(k,j+1)$,
for 
$1 \leq j < J_k$,
\item[(3)]
$\bigcup_{k=1}^{K} Y(k,J_k) = Y$.
\end{itemize}
This implies that for a fixed
$k$,
the union of all
$Y(k,j)$
is invariant under
$\phi$.
It is also clearly closed and so, by minimality, must be all of
$\Sigma$.

\noindent
We shall refer to
$\{Y(k,j) | j = 1, \ldots, J_k\}$
as a tower of height
$J_k$.

\noindent
Now we argue that we can make the partition we have constructed above finer than the 
given one
$\mathcal{P}$,
without changing its essential structure (namely, properties 1-3 above). Suppose
$Z \in \mathcal{P}$
and suppose
$Z$
meets some
$Y(k,j)$
but does not contain it. Divide
$Y(k,j)$
into two clopen sets
$Y(k,j) \cap Z$
and
$Y(k,j) \cap (\Sigma - Z)$.
Unfortunately, this ``disrupts'' the entire
$k$th
tower, so we form
$Y(k,i)^\prime = \phi^{i-j}(Y(k,j) \cap Z)$
and
$Y(k,i)^{\prime\prime} = \phi^{i-j}(Y(k,j) \cap (\Sigma - Z))$,
for each
$i = 1, \ldots, J_k$.
Thus the
$k$th tower breaks into two separate towers (both of height
$J_k$)
with
$Y(k,i)^\prime \subset Z$
and
$Y(k,i)^{\prime\prime}$
disjoint from
$Z$. 
We repeat this for all 
$Z$
and all
$(k,j)$
(which will be a finite process). We then obtain a new 
$K$
and new clopen sets
$Y(k,j)$
(neither will be given a new notation) which satisfy conditions 1-3 above and such that
the partition
$\mathcal{P}^\prime = \{Y(k,j) | k = 1, \ldots, K, j = 1, \ldots J_k\}$
is finer than
$\mathcal{P}$.

\noindent
We are now prepared to define a finite dimensional 
$C^\ast$-
subalgebra of
$C(\Sigma) >\!\!\!\triangleleft_\phi\ \mathbb{Z}$.
In fact, it will be 
$\ast$-
isomorphic to
\[ M_{J_1} \oplus \cdots \oplus M_{J_K}. \]
To do this, it suffices to define matrix units
$e_{ij}^{(k)}$
for all
$k = 1, \ldots, K$
and
$i,j = 1, \ldots, J_k$.
Let
\[ e_{ij}^{(k)} = V^{i-j}\chi_{Y(k,j)} = \chi_{Y(k,i)}V^{i-j}. \]
Using the conditions 1-3, it is routine to check that for fixed
$k$,
$\{e_{ij}^{(k)}\}$
forms a complete system of matrix units for
$M_{J_k}$,
that the projections
\[ p_k = \sum_{i=1}^{J_k} e_{ii}^{(k)} \]
are pairwise orthogonal and sum to the identity. We also note that
$\mathrm{span}\{e_{ii}^{(k)} | k = 1, \ldots, K, i = 1, \ldots, J_k\} = 
 \mathcal{C}(\mathcal{P}^\prime) \supset \mathcal{C}(\mathcal{P})$ 
and that
\[ V\chi_{\Sigma-Y} = \sum_{k=1}^K \sum_{i=2}^{J_k} e_{i,i-1}^{(k)}. \]
The
$C^\ast$-
algebra generated by
$\mathcal{C}(\mathcal{P})$
and
$V\chi_{\Sigma-Y}$
is contained in the finite dimensional algebra we have just described and therefore 
must itself be finite-dimensional. 
\end{proof}

\noindent
We will denote the
$C^\ast$-
algebra generated by
$\mathcal{C}(\mathcal{P})$
and
$V\chi_{\Sigma-Y}$
by
$A(Y,\mathcal{P})$.
\begin{lem} \label{refinePY-lem}
Let
$Y_1$
and
$Y_2$
be two non-empty clopen subsets of
$\Sigma$, 
and let
$\mathcal{P}_1$
and
$\mathcal{P}_2$
be two partitions of
$\Sigma$.
If
$\mathcal{P}_1 \leq \mathcal{P}_2$,
$\chi_{Y_1} \in \mathcal{C}(\mathcal{P}_2)$
and
$Y_1 \supset Y_2$,
then
$A(Y_1, \mathcal{P}_1) \subset A(Y_2, \mathcal{P}_2)$.
\end{lem}
\begin{proof}
Clearly
$\mathcal{C}(\mathcal{P}_1) \subset \mathcal{C}(\mathcal{P}_2)$
and, since
$Y_2 \subset Y_1$,
\[ V\chi_{\Sigma-Y_1} V\chi_{\Sigma-Y_2}\chi_{\Sigma-Y_1} \in A(Y_2, \mathcal{P}_2). \]
\end{proof}

\noindent
\begin{thm}
Let
$Y$
be a non-empty closed subset of
$\Sigma$.
Then
$A_Y$,
the
$C^\ast$-
subalgebra of
$C(\Sigma) >\!\!\!\triangleleft_\phi\ \mathbb{Z}$
generated by
$C(\Sigma)$
and
$VC_0(\Sigma - Y)$,
is an AF-algebra. 
\end{thm}
\begin{proof}
We begin by selecting an increasing sequence of partitions of
$\Sigma$,
$\mathcal{P}_1 \leq \mathcal{P}_2 \leq \cdots$,
whose union generates the topology of
$\Sigma$.
We also choose a decreasing sequence of clopen subsets of
$\Sigma$,
$Y_1 \supset Y_2 \supset \cdots$,
whose intersection is
$Y$.
We will inductively define partitions
$\mathcal{P}_n^\prime$
and finite dimensional subalgebras
$A_n = A(Y_n, \mathcal{P}_n^\prime)$,
for each positive integer
$n$.
Let
$\mathcal{P}_1^\prime = \mathcal{P}_1$
and
$A_1 = A(Y_1, \mathcal{P}_1)$. 
Now assume that we have defined
$\mathcal{P}_n^\prime$
and
$A_n = A(Y_n, \mathcal{P}_n^\prime)$.
We let
$\mathcal{P}_{n+1}^\prime = \mathcal{P}_n^\prime \vee \mathcal{P}_{n+1} \vee 
  \{Y_n, X-Y_n\}$.
Then we have 
$\mathcal{P}_{n+1}^\prime \geq \mathcal{P}_{n+1}$,
$\mathcal{P}_{n+1}^\prime \geq \mathcal{P}_n^\prime$
and
$\chi_{X-Y_n} \in \mathcal{C}(\mathcal{P}_{n+1}^\prime)$.
Let
$A_{n+1} = A(Y_{n+1},\mathcal{P}_{n+1}^\prime)$.

\noindent
We claim that the
$A_n$'s
form a nested sequence of finite dimensional subalgebras of 
$A_Y$
whose union is dense in
$A_Y$. 
First of all,
$\mathcal{C}(\mathcal{P}_n^\prime) \subset C(\Sigma)$
and
$V\chi_{X-Y_n} \in VC_0(\Sigma-Y)$,
since
$Y \subset Y_n$,
so
$A_n \subset A_Y$.
From the properties of
$\mathcal{P}_{n+1}^\prime$
as described in the last paragraph and lemma~\ref{refinePY-lem}, we see that 
$A_n \subset A_{n+1}$,
for all
$n$.
Since the union of the
$\mathcal{P}_n$'s
generates the topology of
$\Sigma$
and
$\mathcal{C}(\mathcal{P}_n) \subset \mathcal{C}(\mathcal{P}_n^\prime) \subset A_n$,
we know that
$C(\Sigma) \subset \overline{\bigcup_n A_n}$.
As
$Y$
is the intersection of the
$Y_n$'s
it is clear that
$VC_0(\Sigma-Y) \subset \overline{\bigcup_n A_n}$.
\end{proof}

\noindent
To get a better understanding of the ordered 
$K_0$-
group of
$C_0(M_0) >\!\!\!\triangleleft_{F_T}\ \mathbb{R}$
Putnam~\cite{Put92} constructed another surface
$\widetilde{M}$
and flow
$\widetilde{F}$
with a proper surjective map
$p: \widetilde{M} \rightarrow M$
which is equivariant for 
$F$
and 
$\widetilde{F}$:
Take a simple closed curve and consider the induced interval exchange transformation
$T$. 
Let
$(M, F^t)$
be the flow constructed from 
$T$
as in section~\ref{IET-sec}. By lemma~\ref{1-form-IET-lem} the flow
$F^t$
is topologically equivalent to the original flow
$f^t$. 

\noindent
Now Putnam splitted the 
$F$-\!
orbit of
$(T(0),\frac{1}{2})$
into 2 parallel orbits. Note that the
$F$-\!
orbit of
$(T(0),\frac{1}{2})$ 
is
\[ \{(T^k(0),s)\ |\ k \geq 1,\  T^k(0) \neq T(0), s \in (0,1]\} - 
   \{(T^l(0),1)\ |\ T^{l+1}(0) = T(0)\}, \]
since as
$t$
approaches
$-\infty$,
$F((T(0),\frac{1}{2}),t)$
converges to
\[ (T(0),0) = (\beta^\prime(\sigma(1)-1),0). \] 
Let
$D_0(T) = \{T^k(0)\ |\ k \geq 1, T^k(0) \neq T(0)\}$.
Let
\[ \Sigma_0 = [0,1] - D_0(T) \cup \{x^+,x^-\ |\ x \in D_0(T)\} \]
with the obvious linear order on
$\Sigma_0$,
using
$x^- < x^+$,
for all
$x \in D_0(T)$.
If
$D_0(T)$ 
is dense,
$\Sigma_0$
endowed with the order topology is a Cantor set.

\noindent
Let
$\widetilde{M}$
be the compact set obtained from
$\Sigma_0 \times [0,1]$
as follows. First, identify the points
$(T(0)^+,0)$
and
$(T(0)^-,0)$
and denote the resulting point by
$(T(0),0)$. 
Identify each of
$\{0\} \times [0,1]$
and
$\{1\} \times [0,1]$
to a point. Finally for each 
$i = 1, \ldots , n$
identify
$[\beta(i-1),\beta(i)] \times \{1\}$
with
$[\beta^\prime(\sigma(i)-1),\beta^\prime(\sigma(i))] \times \{0\}$
via
$T$ --
meaning that for
$x$
in
$D_0$,
we identify
$(x^+,1)$
and
$(x^-,1)$
with
$(T(x)^+,0)$
and
$(T(x)^-,0)$,
respectively. Let
$\widetilde{F}$
be the obvious vertical flow stopped at the images of
$(\beta(i),1)$,
$i = 0,1, \ldots, n$
in
$\widetilde{M}$. 
There is an obvious surjection
$\pi: \widetilde{M} \rightarrow M$,
and we may choose
$\widetilde{F}$
such that
$\pi$
is equivariant, i.e. 
\[ F \circ \textrm{id} \times \pi = \pi \circ \widetilde{F} \]
as maps from
$\mathbb{R} \times \widetilde{M}$
to
$M$.

\noindent
This construction is useful because of
\begin{thm}[{\cite[Thm. 4.1]{Put92}}] \label{K0M0tilde-thm}
Let
$\widetilde{M}_0 = \pi^{-1}(M_0)$. Then
\begin{itemize}
\item[(i)]
$C_0(\widetilde{M}_0)>\!\!\!\triangleleft_{\widetilde{F}}\ \mathbb{R}$
is an AF-algebra, and
\item[(ii)]
$\pi_\ast: K_0(C_0(M_0) >\!\!\!\triangleleft_F\ \mathbb{R}) \rightarrow 
           K_0(C_0(\widetilde{M}_0)>\!\!\!\triangleleft_{\widetilde{F}}\ \mathbb{R})$
is an order isomorphism.
\end{itemize}
\end{thm}

\begin{rem} \label{order-rem}
The order on
$K_0(C_0(M_0) >\!\!\!\triangleleft_F\ \mathbb{R})$
is pulled back from the order on
$K_0(C_0(\widetilde{M}_0)>\!\!\!\triangleleft_{\widetilde{F}}\ \mathbb{R})$
via the statement (contained in the proof) that every projection on
$C_0(\widetilde{M}_0)>\!\!\!\triangleleft_{\widetilde{F}}\ \mathbb{R}$
comes from a projection in
$C_0(M_0) >\!\!\!\triangleleft_F\ \mathbb{R}$.
Since
$C_0(\widetilde{M}_0)>\!\!\!\triangleleft_{\widetilde{F}}\ \mathbb{R}$
is an AF-algebra, the equivalence classes of these projections generate a cone in
$K_0(C_0(\widetilde{M}_0)>\!\!\!\triangleleft_{\widetilde{F}}\ \mathbb{R})$,
hence an order, and similarly for
$K_0(C_0(M_0) >\!\!\!\triangleleft_F\ \mathbb{R})$.
\end{rem}

\noindent
To get a Bratteli diagram of the AF-algebra
$C_0(\widetilde{M}_0)>\!\!\!\triangleleft_{\widetilde{F}}\ \mathbb{R}$
(and hence a description of the ordered
$K_0$-
group)
we need to have a close look at the proof of (the first part of) this theorem: The idea 
is to construct a sequence of open subsets
\[ \widetilde{U}_1 \subset \widetilde{U}_2 \subset \ldots \subset \widetilde{M}_0 \]
exhausting 
$\widetilde{M}_0$
such that the crossed product
$C^\ast$-
algebras 
$C_0(\widetilde{U}_i)>\!\!\!\triangleleft_{\widetilde{F}^{\widetilde{U}_i}}\ \mathbb{R}$
are AF-algebras. Since the union of the
$C_0(\widetilde{U}_i)>\!\!\!\triangleleft_{\widetilde{F}^{\widetilde{U}_i}}\ \mathbb{R}$
is dense in
$C_0(\widetilde{M}_0)>\!\!\!\triangleleft_{\widetilde{F}}\ \mathbb{R}$
we get that this
$C^\ast$-
algebra is AF, too. The set
$\widetilde{U}_1$
is chosen first. Then we describe an iterative procedure for obtaining
$\widetilde{U}_{j+1}$
from
$\widetilde{U}_j$.

\noindent
Choose an integer
$K$
sufficiently large so that, for each
$i = 1, \ldots, n$
\[ \{T^k(0) | 1 \leq k \leq K\} \cap I(i) \]
has at least two points. In the case
$\beta(j-1) < T^K(0) < \beta(j)$
and
$\sigma(j) = n$
and in the case
$\beta(j) < T^K(0) < \beta(j+1)$
and
$\sigma(j+1) = 1$,
replace
$K$
by 
$K+1$.
($T^K(0) = \beta(j)$
is impossible because of the IDOC.)
For 
$i = 0, 1, \ldots, n-1$
let
\[ x(1,i) = \inf(\{T^k(0) | 1 \leq k \leq K\} \cap I(i+1)) \]
and for 
$i = 1, \ldots, n$
let
\[ x(0,i) = \inf(\{T^k(0) | 1 \leq k \leq K\} \cap I(i)). \]
For convenience, set
$\Omega = \{0,1\} \times \{1, \ldots, n-1\} \cup \{(1,0),(0,n)\}$.
For each
$\omega \in \Omega$,
let
$k(\omega)$
be the positive integer such that
$x(\omega) = T^{k(\omega)}(0)$.
We also define
\[ x^\prime(1,i) = \inf(\{T^k(0) | 2 \leq k \leq K+1\} \cap I^\prime(i+1)), \]
\[ x^\prime(0,i) = \inf(\{T^k(0) | 2 \leq k \leq K+1\} \cap I^\prime(i)) \]
for appropriate
$i$,
and
$k^\prime(\omega)$
is given by
$T^{k^\prime(\omega)} = x^\prime(\omega)$
for
$\omega \in \Omega$.
Note that
\[ T\{x(\omega) | \omega \in \Omega\} = \{x^\prime(\omega) | \omega \in \Omega\}. \]
We define closed sets
$Z(i)$,
$Z^\prime(i)$
in
$P = \Sigma_0 \times [0,1]$
by
\[ \begin{array}{l}
   Z(0) = Z^\prime(0) = [0,x(1,0)^-] \times [0,1] \\
   Z(n) = Z^\prime(n) = [x(0,n)^+,1] \times [0,1] \\
   Z(i) = [x(0,i)^+,x(1,i)^-] \times [\frac{3}{4},1] \\
   Z^\prime(i) = [x^\prime(0,i)^+,x^\prime(1,i)^-] \times [0,\frac{1}{4}] 
   \end{array} \]
for
$i = 1, \ldots, n-1$.
We let
$\widetilde{U}$ 
($= \widetilde{U_1}$)
be the complement of the union of the images of the
$Z(i)$,
$Z^\prime(i)$
under the quotient map
$\pi: P \rightarrow \widetilde{M}$.
This open set is 
$\widetilde{M}$
minus the union of closed neighborhoods around the singularities of
$F$.

\noindent
Now, Putnam proved that the
$F^{\widetilde{U}}$-
orbits of the points 
\[ \{(T(0)^\pm,\frac{1}{2}), (x^\prime(\delta,i)^\pm,\frac{1}{2}) | \delta = 0,1; 
     i = 1, \ldots, n-1\} \subset \widetilde{U} \]
are pairwise disjoint. Moreover, these orbits bound 
$n$
strips 
$\cong \Sigma_i \times \mathbb{R}$,
which are
$F^{\widetilde{U}}$-
invariant 
($\Sigma_i$
a Cantor set). These
$F^{\widetilde{U}}$-
flows are topologically conjugated to
$\mathrm{id} \times \tau$
on
$\Sigma_i \times \mathbb{R}$,
where
$\tau$
denotes the canonical flow on
$\mathbb{R}$
(i.e. translation). Consequently,
\[ C_0(\widetilde{U}) >\!\!\!\triangleleft_{F^{\widetilde{U}}} \mathbb{R} \cong
   (\bigoplus_{i=1}^n C_0(\Sigma_i \times \mathbb{R}))  
   >\!\!\!\triangleleft_{\mathrm{id} \times \tau} \mathbb{R} \cong 
   (\bigoplus_{i=1}^n C(\Sigma_i)) \otimes \mathcal{K} \]
(the last isomorphism follows from the approximation of compact operators by 
Hilbert-Schmidt operators on
$L^2(\mathbb{R})$
which arise during the construction of the crossed product).

\noindent
Next, we show how to construct
$\widetilde{U}_2$
satisfying the same conditions: Let
$y$
be one of the points of the set
\[ \{(T(0)^\pm,\frac{1}{2}), (x^\prime(\delta,i)^\pm,\frac{1}{2}) | \delta = 0,1; 
     i = 1, \ldots, n-1\} \]
above. Define 
$y^+ = \lim_{t \rightarrow \infty} F^{\widetilde{U}}(t,y)$. 
Putnam showed that each of the corners of the
$Z(i)$'s 
is a
$y^+$
for some
$y = (x^\prime(\delta,i)^\pm, \frac{1}{2})$
in our set. These
$y$'s
cover all but one pair
$(\delta,i)$,
and Putnam proved that
$y = (x^\prime(\delta,i)^\pm, \frac{1}{2})$
has
$y^+$
contained in the interior of the line segment
$[x(0,j)^+,x(1,j)^-] \times \frac{3}{4}$,
for some
$j$.
In fact, 
$(\delta,i)$
is such that
\[ k^\prime(\delta,i) = \sup\{k^\prime(\delta^\prime,i^\prime) | 
                              \delta^\prime = 0,1, i^\prime = 1, \ldots, n-1\}. \]
Let
$y^+ = (T^l(0),\frac{3}{4})$,
for some
$l \geq 1$. 
Notice that
$T^l(0) = \beta(j)$
is not possible because of the IDOC. Redefine the
$x(\delta,i)$
etc. by replacing
$K$
by
$l$. 
In the case
$T^l(0) < \beta(j)$
and
$\sigma(j) = n$
and in the case
$\beta(j) < T^l(0)$
and
$\sigma(j+1) = 1$,
replace
$K$
by 
$l+1$.
Now, define
$Z(i) = [x(0,i),x(1,i)] \times [\frac{7}{8},1]$
and
$Z^\prime(i) = [x^\prime(0,i),x^\prime(1,i)] \times [0,\frac{1}{8}]$.
It is clear that
$U_2$,
defined as before, has the same properties as 
$U_1$
and
\[ C_0(\widetilde{U}_2) >\!\!\!\triangleleft_{F^{\widetilde{U}_2}} \mathbb{R} \cong
   (\bigoplus_{i=1}^n C_0(\Sigma_i^{(2)})) \otimes \mathcal{K}. \]
It is clear how to continue this process inductively to obtain
$\widetilde{U}_3, \widetilde{U}_4, \ldots $.
From the fact that the orbit of
$0$
is dense in
$[0,\beta(n))$
it follows that the union of all
$\widetilde{U}_j$
will be
$\widetilde{M}_0$.

\noindent
Instead of going into the details of the proof, we present an example for these 
arguments, which hopefully enables the reader to produce a proof by himself and also 
clarifies the meaning of the exceptions
$T^l(0) < \beta(j)$
and
$\sigma(j) = n$
resp.
$\beta(j) < T^l(0)$
and
$\sigma(j+1) = 1$. 
We start with the interval exchange transformation
$T = T(\sigma,\alpha)$
where
$\sigma = (12)$
and
$\alpha = (\sqrt{2}-1,2-\sqrt{2})$. 
In the first step, if we don't include the exceptional cases, we set
$K = 4$
and the second strip does not end in one of the
$Z(i)$'s,
but in
$[x^\prime(0,n)^+,1] \times [0,1]$.
And in the next step, there is no
$y = T^k(0)^\pm$,
$1 \leq k \leq K+1$,
whose orbit ends in the interior of one of the intervals
$[x(0,j)^+,x(1,j)^-] \times \{\frac{3}{4}\}$:

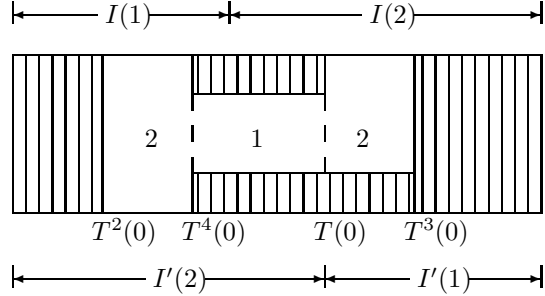
\begin{center}
\begin{picture}(200,110)(0,-30)
\put(0,0){\line(1,0){200}}
\put(0,0){\line(0,1){60}}
\put(0,60){\line(1,0){200}}
\put(200,0){\line(0,1){60}}

\put(0,80){\line(0,-1){8}}
\put(82,80){\line(0,-1){8}}
\put(200,80){\line(0,-1){8}}
\put(32,75){\vector(-1,0){32}}
\put(52,75){\vector(1,0){30}}
\put(35,72){$I(1)$}
\put(132,75){\vector(-1,0){50}}
\put(152,75){\vector(1,0){48}}
\put(135,72){$I(2)$}

\put(0,-20){\line(0,-1){8}}
\put(118,-20){\line(0,-1){8}}
\put(200,-20){\line(0,-1){8}}
\put(50,-25){\vector(-1,0){50}}
\put(73,-25){\vector(1,0){45}}
\put(53,-28){$I^\prime(2)$}
\put(150,-25){\vector(-1,0){32}}
\put(173,-25){\vector(1,0){27}}
\put(153,-28){$I^\prime(1)$}

\put(30,-10){$T^2(0)$}
\put(64,-10){$T^4(0)$}
\put(114,-10){$T(0)$}
\put(148,-10){$T^3(0)$}

\put(50,25){$2$}
\put(90,25){$1$}
\put(130,25){$2$}

\put(34,0){\line(0,1){60}}
\multiput(5,0)(5,0){6}{\line(0,1){60}}

\put(68,45){\line(1,0){50}}
\put(68,45){\line(0,1){15}}
\put(118,45){\line(0,1){15}}
\multiput(70,45)(5,0){10}{\line(0,1){15}}

\put(68,15){\line(1,0){84}}
\put(68,0){\line(0,1){15}}
\put(152,0){\line(0,1){15}}
\multiput(70,0)(5,0){17}{\line(0,1){15}}

\put(152,0){\line(0,1){60}}
\multiput(155,0)(5,0){10}{\line(0,1){60}}

\multiput(68,15)(0,12){3}{\line(0,1){6}}
\multiput(118,15)(0,12){3}{\line(0,1){6}}
\end{picture}
\end{center}

\noindent
Hence we set
$K = 5$ 
and get

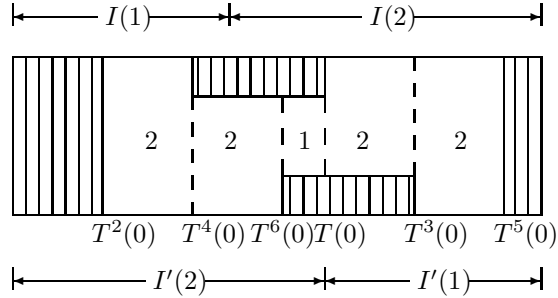
\begin{center}
\begin{picture}(200,110)(0,-30)
\put(0,0){\line(1,0){200}}
\put(0,0){\line(0,1){60}}
\put(0,60){\line(1,0){200}}
\put(200,0){\line(0,1){60}}

\put(0,80){\line(0,-1){8}}
\put(82,80){\line(0,-1){8}}
\put(200,80){\line(0,-1){8}}
\put(32,75){\vector(-1,0){32}}
\put(52,75){\vector(1,0){30}}
\put(35,72){$I(1)$}
\put(132,75){\vector(-1,0){50}}
\put(152,75){\vector(1,0){48}}
\put(135,72){$I(2)$}

\put(0,-20){\line(0,-1){8}}
\put(118,-20){\line(0,-1){8}}
\put(200,-20){\line(0,-1){8}}
\put(50,-25){\vector(-1,0){50}}
\put(73,-25){\vector(1,0){45}}
\put(53,-28){$I^\prime(2)$}
\put(150,-25){\vector(-1,0){32}}
\put(173,-25){\vector(1,0){27}}
\put(153,-28){$I^\prime(1)$}

\put(30,-10){$T^2(0)$}
\put(64,-10){$T^4(0)$}
\put(90,-10){$T^6(0)$}
\put(114,-10){$T(0)$}
\put(148,-10){$T^3(0)$}
\put(182,-10){$T^5(0)$}

\put(50,25){$2$}
\put(80,25){$2$}
\put(108,25){$1$}
\put(130,25){$2$}
\put(167,25){$2$}

\put(34,0){\line(0,1){60}}
\multiput(5,0)(5,0){6}{\line(0,1){60}}

\put(68,45){\line(1,0){50}}
\put(68,45){\line(0,1){15}}
\put(118,45){\line(0,1){15}}
\multiput(70,45)(5,0){10}{\line(0,1){15}}

\put(102,15){\line(1,0){50}}
\put(102,0){\line(0,1){15}}
\put(152,0){\line(0,1){15}}
\multiput(105,0)(5,0){10}{\line(0,1){15}}

\put(186,0){\line(0,1){60}}
\multiput(190,0)(5,0){2}{\line(0,1){60}}

\multiput(68,0)(0,10){5}{\line(0,1){5}}
\multiput(102,15)(0,12){3}{\line(0,1){6}}
\multiput(118,15)(0,12){3}{\line(0,1){6}}
\multiput(152,15)(0,10){5}{\line(0,1){5}}
\end{picture}
\end{center}

\noindent
The next step leads to
$K = 7$
and

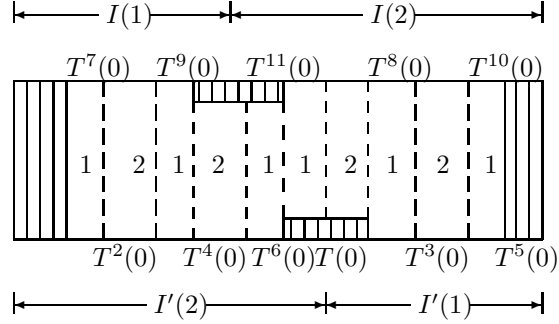
\begin{center}
\begin{picture}(200,120)(0,-30)
\put(0,0){\line(1,0){200}}
\put(0,0){\line(0,1){60}}
\put(0,60){\line(1,0){200}}
\put(200,0){\line(0,1){60}}

\put(0,90){\line(0,-1){8}}
\put(82,90){\line(0,-1){8}}
\put(200,90){\line(0,-1){8}}
\put(32,85){\vector(-1,0){32}}
\put(52,85){\vector(1,0){30}}
\put(35,82){$I(1)$}
\put(132,85){\vector(-1,0){50}}
\put(152,85){\vector(1,0){48}}
\put(135,82){$I(2)$}

\put(0,-20){\line(0,-1){8}}
\put(118,-20){\line(0,-1){8}}
\put(200,-20){\line(0,-1){8}}
\put(50,-25){\vector(-1,0){50}}
\put(73,-25){\vector(1,0){45}}
\put(53,-28){$I^\prime(2)$}
\put(150,-25){\vector(-1,0){32}}
\put(173,-25){\vector(1,0){27}}
\put(153,-28){$I^\prime(1)$}

\put(20,62){$T^7(0)$}
\put(30,-10){$T^2(0)$}
\put(54,62){$T^9(0)$}
\put(64,-10){$T^4(0)$}
\put(88,62){$T^{11}(0)$}
\put(90,-10){$T^6(0)$}
\put(114,-10){$T(0)$}
\put(134,62){$T^8(0)$}
\put(148,-10){$T^3(0)$}
\put(172,62){$T^{10}(0)$}
\put(182,-10){$T^5(0)$}

\put(25,25){$1$}
\put(45,25){$2$}
\put(60,25){$1$}
\put(75,25){$2$}
\put(94,25){$1$}
\put(108,25){$1$}
\put(125,25){$2$}
\put(141,25){$1$}
\put(160,25){$2$}
\put(178,25){$1$}

\put(20,0){\line(0,1){60}}
\multiput(5,0)(5,0){4}{\line(0,1){60}}

\put(68,52){\line(1,0){34}}
\put(68,52){\line(0,1){8}}
\put(102,52){\line(0,1){8}}
\multiput(70,52)(5,0){7}{\line(0,1){8}}

\put(102,8){\line(1,0){32}}
\put(102,0){\line(0,1){8}}
\put(134,0){\line(0,1){8}}
\multiput(105,0)(5,0){6}{\line(0,1){8}}

\put(186,0){\line(0,1){60}}
\multiput(190,0)(5,0){2}{\line(0,1){60}}

\multiput(34,0)(0,9){7}{\line(0,1){6}}
\multiput(54,0)(0,9){7}{\line(0,1){6}}
\multiput(68,0)(0,8){7}{\line(0,1){4}}
\multiput(88,0)(0,8){7}{\line(0,1){4}}
\multiput(102,8)(0,10){5}{\line(0,1){5}}
\multiput(118,8)(0,8){7}{\line(0,1){4}}
\multiput(134,8)(0,8){7}{\line(0,1){4}}
\multiput(152,0)(0,9){7}{\line(0,1){6}}
\multiput(172,0)(0,9){7}{\line(0,1){6}}
\end{picture}
\end{center}

\noindent
The whole construction tells us that
$C_0(\widetilde{M}_0)>\!\!\!\triangleleft_{\widetilde{F}}\ \mathbb{R}$
is the inductive limit of 
$C^\ast$-
algebras of type
$(\bigoplus_{i=1}^n C(\Sigma_i)) \otimes \mathcal{K}$.
But these 
$C^\ast$-
algebras are already AF-algebras because the
$\Sigma_i$
are Cantor sets (\cite[Ex.III.2.5]{Dav96}). Hence,
$C_0(\widetilde{M}_0)>\!\!\!\triangleleft_{\widetilde{F}}\ \mathbb{R}$
is also an AF-algebra. It can be presented as the inductive limit of 
finite-dimensional
$C^\ast$-
algebras in the following way: Let
$A_j$
be the 
$C^\ast$-
subalgebra of
$C_0(\widetilde{U}_j) >\!\!\!\triangleleft_{\widetilde{F}^j} \mathbb{R}$
which is mapped to
$\bigoplus_{i=1}^n \mathbb{C} \otimes \mathcal{K}$ 
under the isomorphism above, for
$j = 1, 2, \ldots$. 
It is not hard to check that
$A_j \subset A_{j+1}$
for all
$j$
and that the union of the
$A_j$
is dense in 
$C_0(\widetilde{M}_0)>\!\!\!\triangleleft_{\widetilde{F}}\ \mathbb{R}$.
Obviously,
$K_0(A_j) \cong \mathbb{Z}^n$,
for
all
$j$,
and the embeddings induce isomorphisms of abelian groups given by a matrix of 
$1$'s
in the diagonal, one off-diagonal
$1$
and the rest of the entries 
$0$.
This form of the matrix reflects the fact that in every step of the construction above,
exactly one strip is divided in two strips and one of these two strips becomes part of 
an already existing one.

\section{Order-isomorphisms of $K_0$-groups} \label{isoK0-sec}

\noindent
The aim is now to show that the
$K_0$-
groups of
$C(\Sigma) >\!\!\!\triangleleft_{\phi_F}\ \mathbb{Z}$
and the crossed product 
$C^\ast$-\!
algebras
$C_0(M_0) >\!\!\!\triangleleft_F\ \mathbb{R}$
and
$C(M) >\!\!\!\triangleleft_F\ \mathbb{R}$
are order-isomorphic. 

\noindent
As a first step, it is an easy exercise in K-theory to prove
\begin{prop} \label{K0Flow-prop}
The 
$K_0$-
groups of
$C_0(M_0) >\!\!\!\triangleleft_F\ \mathbb{R}$
and
$C(M) >\!\!\!\triangleleft_F\ \mathbb{R}$ 
are order isomorphic.
\end{prop}
\begin{proof}
If
$\{x_1, \ldots, x_l\}$
are the singularities of the minimal flow
$F$,
we have the short exact sequence
\[ 0 \rightarrow C_0(M_0) >\!\!\!\triangleleft_F\ \mathbb{R} \stackrel{\phi}{\rightarrow}
                 C(M) >\!\!\!\triangleleft_F\ \mathbb{R} \stackrel{\psi}{\rightarrow} 
                 \bigoplus_{i=1}^l C_0(\mathbb{R}) \rightarrow 0. \]
The long exact sequence of K-functors yields
\[ \begin{array}{rcccl}
   K_2(\bigoplus_{i=1}^l C_0(\mathbb{R})) & \rightarrow & 
   K_1(C_0(M_0) >\!\!\!\triangleleft_F\ \mathbb{R}) & \rightarrow &
   K_1(C(M) >\!\!\!\triangleleft_F\ \mathbb{R}) \rightarrow \\
    & & \\
    & \rightarrow & K_1(\bigoplus_{i=1}^l C_0(\mathbb{R})) & \rightarrow & 
                    K_0(C_0(M_0) >\!\!\!\triangleleft_F\ \mathbb{R}) \rightarrow \\
    & & \\
    & \rightarrow & K_0(C(M) >\!\!\!\triangleleft_F\ \mathbb{R}) & \rightarrow &
                    K_0(\bigoplus_{i=1}^l C_0(\mathbb{R})). 
   \end{array}\]
Next,
$C_0(\mathbb{R}) = C_0(\mathbb{R}) \otimes \mathbb{C}$
is the suspension
$\Sigma\mathbb{C}$
of the
$C^\ast$-\!
algebra
$\mathbb{C}$
\cite[Def.1.15]{GVF-NC}.
By using
$K_{n+1}(A) = K_n(\Sigma A)$
\cite[Prop.3.26]{GVF-NC},
$K_1(\mathbb{C}) = 0$
\cite[p.128]{GVF-NC} and the split exact sequences associated to the direct sum
$\bigoplus_{i=1}^l K_0(C_0(\mathbb{R}))$
\cite[Prop.3.29]{GVF-NC}, we get
\[ \begin{array}{rclcl}
   K_0(\bigoplus_{i=1}^l C_0(\mathbb{R})) & = & 0, & & \\ 
    & & \\
   K_1(\bigoplus_{i=1}^l C_0(\mathbb{R})) & =  & \bigoplus_{i=1}^l K_1(C_0(\mathbb{R})) 
    & = & 
   \bigoplus_{i=1}^l K_1(\Sigma\mathbb{C}) = \bigoplus_{i=1}^l K_2(\mathbb{C}) = \\ 
    & & \\
    & = & \bigoplus_{i=1}^l K_0(\mathbb{C}) & =  & \mathbb{Z}^l\ \ \mathrm{and} \\
    & & \\ 
   K_2(\bigoplus_{i=1}^l C_0(\mathbb{R})) & =  & K_0(\bigoplus_{i=1}^l C_0(\mathbb{R}) 
   & = & 0, 
   \end{array}\]
by Bott periodicity \cite[Thm.3.34]{GVF-NC} and 
$K_0(\mathbb{C}) = \mathbb{Z}$
\cite[p.96]{GVF-NC}.
Connes' analogue to the Thom isomorphism \cite[II.C, Thm.8]{NCG} tells us that
\[ \begin{array}{l}
   K_1(C_0(M_0) >\!\!\!\triangleleft_F\ \mathbb{R}) \cong K_2(C_0(M_0)) \cong 
   K_0(C_0(M_0)) \\
     \\
   K_1(C(M) >\!\!\!\triangleleft_F\ \mathbb{R}) \cong K_2(C(M)) \cong K_0(C(M)). 
   \end{array} \]
Since the second map in the short exact sequence
\[ 0 \rightarrow C_0(M_0) \rightarrow C(M) \rightarrow C(\{x_1, \ldots, x_l\}) 
      \rightarrow 0 \]
has a (in fact, many) splitting(s), we have again by \cite[Prop.3.29]{GVF-NC} that there 
is a short exact sequence
\[ 0 \rightarrow K_0(C_0(M_0)) \rightarrow K_0(C(M)) \rightarrow 
                 K_0(C(\{x_1, \ldots, x_l\}) = \bigoplus_{i=1}^l K_0(\mathbb{C}) 
      \rightarrow 0. \]
Consequently, 
\[  \phi_\ast: K_0(C_0(M_0) >\!\!\!\triangleleft_F\ \mathbb{R}) 
    \stackrel{\cong}{\rightarrow}
               K_0(C(M) >\!\!\!\triangleleft_F\ \mathbb{R}). \]
is an isomorphism. 

\noindent
We need some additional arguments to show that the 
$K_0-$
groups of the considered
$C^\ast$-
algebras have an order. As described in remark~\ref{order-rem} the order on
$K_0(C_0(M_0) >\!\!\!\!\!\triangleleft_F\ \mathbb{R})$
is pulled back from the order on
$K_0(C_0(\widetilde{M}_0)>\!\!\!\!\!\triangleleft_{\widetilde{F}}\ \mathbb{R})$. 
It remains to construct a projection
$p \in \mathcal{P}(C_0(M_0) >\!\!\!\!\!\triangleleft_F\ \mathbb{R})$
for every projection
$q \in \mathcal{P}(C(M) >\!\!\!\triangleleft_F\ \mathbb{R})$
such that
\[ \phi_\ast([p]) = [q]. \]
But this is trivial, because the image
$\psi(q)$
is a projection in some matrix algebra over
$\bigoplus_{i=1}^l C_0(\mathbb{R})$,
hence equivalent to
$0$,
and consequently,
\[ q \in C_0(M_0) >\!\!\!\triangleleft_F\ \mathbb{R}. \]
\end{proof}

\noindent
Next we prove
\begin{thm} \label{FlowCantor-thm}
Let
$(\Sigma,\phi_F)$
be the dynamical system on a Cantor set
$\Sigma$
corresponding to the minimal flow 
$F$
induced by (the imaginary part of) a holomorphic
$1$-
form on a Riemann surface
$X = M$.
Then 
\[ (K_0(C_0(M_0) >\!\!\!\triangleleft_F\ \mathbb{R}), 
    K_0^+(C_0(M_0) >\!\!\!\triangleleft_F\ \mathbb{R})) \cong 
   (K_0(C(\Sigma) >\!\!\!\triangleleft_{\phi_F}\ \mathbb{Z}), 
    K_0^+(C(\Sigma) >\!\!\!\triangleleft_{\phi_F}\ \mathbb{Z}). \]
\end{thm}
\begin{proof}
The idea is to bring the apparent similarity of the ``method of towers'' and the 
construction for
$K_0(C_0(M_0) >\!\!\!\triangleleft_F\ \mathbb{R})$
in section~\ref{crpr-sec}
into a mathematically exact form. This is done by constructing 
AF-subalgebras
$B_j \subset C(\Sigma) >\!\!\!\triangleleft_{\phi_F}\ \mathbb{Z}$
which are stably isomorphic to
$A_j$,
such that
$B := \lim_{j \rightarrow \infty} B_j \subset 
      C(\Sigma) >\!\!\!\triangleleft_{\phi_F}\ \mathbb{Z}$
is stably isomorphic to 
$\lim_{j \rightarrow \infty} A_j = 
 C_0(\widetilde{M}_0)>\!\!\!\triangleleft_{\widetilde{F}}\ \mathbb{R}$.
Furthermore we show that the inclusion
$i: C(\Sigma) \hookrightarrow C(\Sigma) >\!\!\!\triangleleft_{\phi_F}\ \mathbb{Z}$
factors through 
$B$:

\begin{center}
\mbox{
\xymatrix{
                & C(\Sigma) \ar[ld]_{i_1} \ar[rd]^{i} & \\
B \ar[rr]_{i_2} & & C(\Sigma) >\!\!\!\triangleleft_{\phi_F}\ \mathbb{Z}
         }}
\end{center}

\noindent
Then by \cite[Cor.2.4]{Put90}, for every 
$a$
in
$K_0(C(\Sigma) >\!\!\!\triangleleft_{\phi_F}\ \mathbb{Z})^+$
there is a 
$c$
in
$K_0(C(\Sigma))^+$
such that
$i_\ast(c) = a$.
Letting
$b = i_{1,\ast}(c)$
we have
$i_{2,\ast}(b) = a$.
Hence 
$i_{2,\ast}$
is surjective. It is also injective because
\[ K_0(B) = K_0(C_0(\widetilde{M}_0)>\!\!\!\triangleleft_{\widetilde{F}}\ \mathbb{R}) =
   \mathbb{Z}^n = K_0(C(\Sigma) >\!\!\!\triangleleft_{\phi_F}\ \mathbb{Z}), \]
where the first isomorphism follows from
$B$
stably isomorphic to
$\lim_{j \rightarrow \infty} A_j = 
 C_0(\widetilde{M}_0)>\!\!\!\triangleleft_{\widetilde{F}}\ \mathbb{R}$. 
The conclusion of the theorem follows.

\noindent
The towers leading to the
$B_j$'s 
consist of the intervals
$[T^l(0)^+, T^m(0)^-]$
such that the lines
$T^l(0)^+ \times [0,1]$
and
$T^m(0)^- \times [0,1]$
border the same strip in the
$j$th
step of the construction above. In particular, there is a minimal pair
$(l,m)$
such that the intervals
\[ \begin{array}{rcl}
   Y(k,0) & = & [T^l(0)^+, T^m(0)^-], \\
   Y(k,1) & = & [T^{l+1}(0)^+, T^{m+1}(0)^-] , \\
          & \vdots & \\
   Y(k,J_k) & = & [T^{l+J_k}(0)^+, T^{m+J_k}(0)^-]
   \end{array} \]
form the  
$k$th 
strip, 
$k = 1, \ldots , n$. 
In this picture, the strips
$[0, T^{k(1,0)}(0)^-] \times [0,1]$
and
$[T^{k(0,n)}(0)^+, 1] \times [0,1]$
are missing, but they can be naturally included as the strips preceding those with left 
boundary
$T(0)^+ \times [0,1]$
resp. right boundary
$T(0)^- \times [0,1]$.
We get a partition of
$\Sigma$
into the intervals
$Y(k,j)$,
$k = 1, \ldots , n$,
$j = 0, \ldots J_k$.

\noindent
This partition may be used as in the ``method of towers'' to produce a finite 
dimensional
$C^\ast$-
subalgebra of 
$C(\Sigma) >\!\!\!\triangleleft_{\phi_F}\ \mathbb{Z}$,
\[ M_{J_1} \oplus \ldots \oplus M_{J_n}, \]
defined with matrix units
\[ e_{ij}^{(k)} = V^{i-j}\chi_{Y(k,j)} = \chi_{Y(k,i)}V^{i-j} \]
for all
$k = 1, \ldots , n$,
$i,j = 0, \ldots J_k$.
Let
$B_j$
be this finite 
dimensional
$C^\ast$-
algebra.

\noindent
Since the
$(j+1)$st
partition is finer than the 
$j$th
partition and the union of all partitions generates the toplogy of
$\Sigma$,
the 
$B_j$'s
form a nested sequence of finite 
dimensional
$C^\ast$-
algebras whose inductive limit
$B$
contains
$C(\Sigma)$.

\noindent
On the other hand, we can also carryover the argument that yields the form of the matrix
describing the homomorphism 
$K_0(A_j) \rightarrow K_0(A_{j+1})$
in the natural
$\mathbb{Z}^n$-
bases. Hence,
$B$
and
$\lim_{j \rightarrow \infty} A_j$
are stably isomorphic, and the proof is finished. 
\end{proof}

\noindent
\begin{rem}
It would be nice to have a
$\ast$-
homomorphism between
$C_0(M_0) >\!\!\!\triangleleft_F\ \mathbb{R}$
and
$C(\Sigma) >\!\!\!\triangleleft_{\phi_F}\ \mathbb{Z}$
which induces an order-isomorphism of their 
$K_0$-
groups.
Indeed, Putnam~\cite[Thm.3.2]{Put92} proved the following
\begin{thm}
Let 
$C(\Sigma) >\!\!\!\triangleleft_{\phi_T}\ \mathbb{Z}$
and
$C_0(M_0) >\!\!\!\triangleleft_{F_T}\ \mathbb{R}$
arise from a transformation
$T(\sigma,\alpha)$
exchanging
$n$
intervals satisfying the I.D.O.C. as above. Then there are isomorphisms
\[ K_0(C(\Sigma) >\!\!\!\triangleleft_{\phi_T}\ \mathbb{Z}) \cong \mathbb{Z}^n,\ 
   K_0(C_0(M_0) >\!\!\!\triangleleft_{F_T}\ \mathbb{R}) \cong \mathbb{Z}^n \]
and an injective 
$\ast$-
homomorphism
\[ \rho: C_0(M_0) >\!\!\!\triangleleft_{F_T}\ \mathbb{R} \rightarrow
         (C(\Sigma) >\!\!\!\triangleleft_{\phi_T}\ \mathbb{Z}) \otimes \mathcal{K}, \]
where
$\mathcal{K}$
denotes the 
$C^\ast$-
algebra of compact operators on a separable, infinite-dimensional Hilbert space.

\noindent
Furthermore, using the isomorphisms above,
\[ \rho_\ast: K_0(C_0(M_0) >\!\!\!\triangleleft_{F_T}\ \mathbb{R}) \rightarrow
              K_0(C(\Sigma) >\!\!\!\triangleleft_{\phi_T}\ \mathbb{Z}) \]
is given by multiplication by
$L^\sigma$,
an
$n \times n$
matrix with
\[ L^\sigma_{ij} = \left\{ \begin{array}{cl}
                           1 & \mathrm{if\ } i > j \mathrm{\ and\ } 
                               \sigma(i) < \sigma(j) \\
                           -1 & \mathrm{if\ } i < j \mathrm{\ and\ } 
                               \sigma(i) > \sigma(j) \\
                           0 & \mathrm{otherwise}.
                           \end{array}  \right.\]
\end{thm}

\noindent
Unfortunately, 
$L^\sigma$
need not be invertible (e.g., if
$n$
is odd, since then
$L^\sigma$
describes an alternating quadratic form).

\noindent
It is not known to the author whether there is another
$\ast$-
homomorphism inducing an order-isomorphism between the 
$K_0$-
groups, or whether the two crossed product
$C^\ast$-
algebras
$C_0(M_0) >\!\!\!\triangleleft_F\ \mathbb{R}$
and
$C(\Sigma) >\!\!\!\triangleleft_{\phi_F}\ \mathbb{Z}$
are (strongly) Morita equivalent.
\end{rem}

\noindent
Proposition~\ref{K0Flow-prop} and theorem~\ref{FlowCantor-thm} show that the
$C^\ast$
algebras associated in section~\ref{crpr-sec} to a minimal flow induced by the imaginary 
part of a holomorphic
$1$-
form have the same ordered 
$K_0$-
group. We finally prove that the cone defining the order for this abelian group is dual to
the cone 
$\Sigma(\sigma, \alpha)$
of 
$T$-
invariant Borel measures on
$[0,\beta(n))$
where
$T$
is an induced interval exchange transformation on the interval
$[0,\beta(n))$. 

\noindent
Recall the description of
$\Sigma(\sigma, \alpha)$ 
in proposition~\ref{ergcone-prop}: Let
$J_1 \supset J_2 \supset \cdots$
be a sequence of admissible intervals shrinking to a point. Then
\[ \Sigma(\sigma, \alpha) \cong \bigcap_{i=1}^\infty (A_1 A_2 \cdots A_i \Lambda_n) \]
where the
$A_i$ 
are the matrices describing the transition from the
$(i-1)$st
to the
$i$th 
induced interval exchange transformation and 
$\Lambda_n \subset \mathbb{R}^n$
is the simplex generated by the unit vectors in
$\mathbb{R}^n$.

\noindent
For every point
$y_0 \in [0,\beta_\alpha(n))$
it is possible to choose a sequence of admissible 
subintervals 
$J_k = [a_k,b_k)$
such that
\[ J_1 = [0,\beta_\alpha(n)) \supset J_2 \supset \cdots \{y_0\} \]
contract to this point
$\{y_0\}$: 
Start with the admissible interval
$J_2 = [\beta_\alpha(j), \beta_\alpha(j+1))$
containing
$y_0$.
Take as
$J_3$
the interval
$[\beta_{\alpha^{(2)}}(k), \beta_{\alpha^{(2)}}(k+1))$
of the induced interval exchange transformation 
$\alpha^{(2)}$
containing
$y_0$, 
and iterate this procedure. The
$J_i$
contract to
$y_0$:
Assume that another point
$y_1 \neq y_0$
lies in all the
$J_n$.
Let
$T^{(n)}_{|J_n} = T^{l^{(n)}}$.
Then by construction,
$T^m(y_0)$
and
$T^m(y_1)$
would lie in the same interval of
$\alpha = \alpha_1$,
for all
$0 \leq m \leq l^{(n)}$.
Now,
$l^{(n)}$
increases for increasing
$n$. 
This means that
$T^m(y_0)$
and
$T^m(y_1)$
would lie in the same interval of
$\alpha = \alpha_1$,
for all
$m \in \mathbb{N}$.
But this is impossible: By the minimality assumption, there are 
$i$,
$l \in \mathbb{N}$
such that
$T^{-l}\beta_\alpha(i)$
separates
$y_0$
and
$y_1$.
This certainly implies that
$T^l(y_0)$
and
$T^l(y_1)$
lie in different intervals of
$\alpha$.

\noindent
The sequence of
$\mathbb{Z}$-\!
invertible matrices
$A^{(k)} \in M(n, \mathbb{Z})$
describing the transition from
$\alpha^{(k)}$
to
$\alpha^{(k+1)}$
(that is,
$\alpha^{(k)} = A^{(k)} \alpha^{(k+1)}$)
allow us to define an ordered group
$(G_T, G_T^+)$
as the direct limit of 
\[ (\mathbb{Z}^n,(\mathbb{Z}^n)^+) \stackrel{A^{(1)}}{\rightarrow} 
   (\mathbb{Z}^n,(\mathbb{Z}^n)^+) \stackrel{A^{(2)}}{\rightarrow}
   (\mathbb{Z}^n,(\mathbb{Z}^n)^+) \stackrel{A^{(3)}}{\rightarrow} \cdots . \]
The homomorphisms
$A^{(k)}: (\mathbb{Z}^n,(\mathbb{Z}^n)^+) \rightarrow (\mathbb{Z}^n,(\mathbb{Z}^n)^+)$
are positive, since all the entries of
$A^{(k)}$
are positive integers. Furthermore, this direct limit is dual to the cone
$\bigcap_{i=1}^\infty (A_1 A_2 \cdots A_i \Lambda_n)$.
\begin{prop} \label{CantorK0-prop}
Let
$T(\sigma,\alpha)$
and
$(\Sigma,\phi)$
be the interval exchange transformation resp. dynamical system on a Cantor set
$\Sigma$
corresponding to the minimal flow induced by (the imaginary part of) a holomorphic
$1$-
form on a Riemann surface
$X$.
Then
\[ (G_T, G_T^+) \cong (K_0(C(\Sigma) >\!\!\!\triangleleft_\phi\ \mathbb{Z}), 
                       K_0^+(C(\Sigma) >\!\!\!\triangleleft_\phi\ \mathbb{Z})) \]
as ordered groups.
\end{prop}
\begin{proof}
By theorem~\ref{PutK0-thm} we have to prove that
\[ (G_T, G_T^+) \cong (K_0(A_{\{y\}}),K_0(A_{\{y\}})^+) \]
for some point
$y \in \Sigma$.
It is quite self-evident to set
$y = y_0 \in [0,\beta_\alpha(n)) \subset \Sigma$
and to try to identify the ordered group
$(\mathbb{Z}^n,(\mathbb{Z}^n)^+)$
from the direct limit defining 
$(G_T, G_T^+)$
with the ordered
$K_0$-
group of
$A(Y_k,\mathcal{P}_k)$
for some clopen interval
$Y_k \subset \Sigma$
and partition
$\mathcal{P}_k$
which are in some way derived from the admissible interval
$J_k = [a_k,b_k)$
and the interval exchange transformation
$T(\sigma^{(k)},\alpha^{(k)})$.

\noindent
Again, it is eye-catching to set
$Y_k = [a_k^+, (a_k + \beta_{\alpha^{(k)}}(n))^-] \subset \Sigma$.
This gives a generic choice for the partitions
$\mathcal{P}_k$: 
Define
$\lambda_k : Y_k \rightarrow \mathbb{Z}$
and
$Y_k(l,j)$
as above such that the clopen intervals
$Y_k(l,J_l)$
are the same as the intervals
$[(a_k + \beta_{\alpha^{(k)}}(i-1))^+, (a_k + \beta_{\alpha^{(k)}}(i))^-]$,
$i = 1, \ldots, n$. 
(This means in particular that some of the
$J_l$
may be equal.) Then let
$\mathcal{P}_k$
be the partition of
$\Sigma$
into the towers
$\{Y_k(l,j) | j = 1, \ldots, J_l\}$
of height
$J_l$,
$l = 1, \ldots, n$.
Obviously,
\[ \{Y_k, \Sigma - Y_k\} \leq \mathcal{P}_k \leq \mathcal{P}_{k+1}. \]
Comparing with the arguments above we see that
\[ A_k = A(Y_k, \mathcal{P}_k) \subset M_{J_1} \oplus \cdots \oplus M_{J_n} \]
is generated by the elements (see notation above)
\[ e_{ii}^{(l)},\ l = 1, \ldots, n,\ i = 1, \ldots, J_l;\ \ 
   f = \sum_{l=1}^n \sum_{i=2}^{J_l} e_{ii-1}^{(l)}. \]
Since 
$e_{ii}^{(l)} \cdot f = e_{ii-1}^{(l)}$,
$f^m = \sum_{l=1}^n \sum_{i=m+1}^{J_l} e_{ii-m}^{(l)}$
and
$(e_{ij}^{(l)})^\ast = e_{ji}^{(l)}$,
we have
\[ A_k \cong M_{J_1} \oplus \cdots \oplus M_{J_n}. \]

\noindent
We still have to show that the union of the partitions
$\mathcal{P}_k$
generates the topology of
$\Sigma$.
But this is clear since the lengths of the intervals
$Y_k$
tend to
$0$,
and the intervals
$Y_k(l,j)$
are subsets of
$Y_k$.
Hence
\[ A_{\{y_0\}} = \overline{\bigcup_k A_k}. \]
To compute
$(K_0(A_{\{y_0\}}), K_0^+(A_{\{y_0\}}))$
we need to determine the inclusions
$A_k \rightarrow A_{k+1}$.
It is enough to look at the
$\ast$-
homomorphisms
$M_{J_l^{(k)}} \rightarrow M_{J_m^{(k+1)}}$,
$1 \leq l,m \leq n$,
which are unitarily equivalent to direct sums of the identity map on 
$M_{J_l^{(k)}}$,
(by \cite[III.2.1]{Dav96}). Let
$e_{ij}^{(l,k)}, e_{ij}^{(m,k+1)}$
be the matrix units of
$A_k$
resp.
$A_{k+1}$
from above. The number of identity representations in the
$\ast$-
homomorphisms
$M_{J_l^{(k)}} \rightarrow M_{J_m^{(k+1)}}$
is given by the rank of the image matrix in
$M_{J_m^{(k+1)}}$
of
$e_{J_lJ_l}^{(l,k)} = \chi_{Y_k(l,J_l)}$. 
But the interpretation of the 
$A_{lm}^{(k)}$
in proposition~\ref{indIET-prop} shows that
$Y_k(l,J_l)$
is composed of 
$A_{lm}^{(k)}$
elements of the tower
$\{Y_{k+1}(m,j) | j = 1, \ldots, J_m^{(k+1)}\}$,
$1 \leq m \leq n$,
hence
$A_{lm}^{(k)}$
is the rank of the image matrix of
$e_{J_lJ_l}^{(l,k)}$
in
$M_{J_m^{(k+1)}}$.
In particular, the matrix
$A^{(k)}$
describes the induced map of
$K_0$-
groups.
 
\noindent
We still have to show that the union of the partitions
$\mathcal{P}_k$
generates the topology of
$\Sigma$.
But this is clear since the lengths of the intervals
$Y_k$
tend to
$0$,
and the intervals
$Y_k(l,j)$
are subsets of
$Y_k$.
\end{proof}

\end{document}